\documentclass[12pt]{amsart}
\usepackage{amssymb,amscd}
\usepackage[mathscr]{eucal}
\input epsf
\newtheorem{thm}{Theorem}[section]
\newtheorem{lemma}[thm]{Lemma}
\newtheorem{defin}[thm]{Definition}
\newtheorem{conjecture}[thm]{Conjecture}
\newtheorem{prop}[thm]{Proposition}
\newtheorem{rem}[thm]{Remark}
\newtheorem{cor}[thm]{Corollary}
\begin{document}

\newcommand{\C}{{\Bbb C}}
\newcommand{\HH}{{\Bbb H}}
\newcommand{\Q}{{\Bbb Q}}
\newcommand{\M}{{\mathcal M}}
\newcommand{\R}{{\Bbb R}}
\newcommand{\Z}{{\Bbb Z}}
\newcommand{\E}{{\mathcal E}}
\newcommand{\SU}{\o{SU}(2)}
\newcommand{\SO}{\o{SO}(3)}
\newcommand{\Hom}{\o{Hom}}
\newcommand{\Spin}{\o{Spin}(2)}
\newcommand{\Spinm}{\o{Spin}_-(2)}
\newcommand{\Pin}{\o{Pin}(2)}
\renewcommand{\O}{{\mathcal O}}
\renewcommand{\P}{{\Bbb P}}
\renewcommand{\o}{\operatorname}
\newcommand{\tr}{\o{tr}}
\title
{Topological Dynamics on Moduli Spaces II}
\author{Joseph P. Previte \and Eugene Z. Xia}
\address{ School of Science,
Penn State Erie, The Behrend College,
Erie, PA 16563}
\address{Department of Mathematics \& Statistics,
University of Massachusetts, Amherst, MA 01003-4515}
\email{jpp@vortex.bd.psu.edu {\it (Previte)}, xia@math.umass.edu {\it (Xia)}}
\date{\today} 
\maketitle
\begin{abstract}
Let $M$ be a Riemann surface with boundary  
$\partial M$ and genus greater than zero.  Let
$\Gamma$ be the mapping class group of $M$ fixing $\partial M$.
The group $\Gamma$ acts on 
${\mathcal M}_{\mathcal C} = \Hom_{\mathcal C}(\pi_1(M),\SU)/\SU$ which
is the space of $\SU$-gauge equivalence classes of flat 
$\SU$-connections on $M$ with fixed holonomy on $\partial M$.
We study the topological dynamics of the $\Gamma$-action and 
give conditions for the individual $\Gamma$-orbits to be dense
in ${\mathcal M}_{\mathcal C}$.  
\end{abstract}

\section{Introduction}
Let $M$ be a Riemann surface of genus $g$ with $n$
boundary components (circles).  Let 
$$
\{C_1, C_2,..., C_n\} \subset \pi_1(M)
$$
be elements in the fundamental group that correspond to
these $n$ boundary components.

The space of $\SU$-gauge equivalence classes of 
$\SU$-connections $YM_2(\SU)$ is the well known Yang-Mills two space 
of quantum field theory.  Inside $YM_2(\SU)$ is the moduli space 
${\mathcal M}$ of flat $\SU$-connections.
The moduli space ${\mathcal M}$ has an interpretation that relates
to the representation space $\Hom(\pi_1(M),\SU),$ which is
an algebraic variety.  The group $\SU$ acts on $\Hom(\pi_1(M),\SU)$
by conjugation, and the resulting quotient space is precisely
$$
{\mathcal M} = \Hom(\pi_1(M),\SU)/\SU.  
$$

Note that a conjugacy class in $\SU$ is determined by its trace.
Hence specifying a conjugacy class is the same as specifying
a real number in $[-2,2]$. 
To each $C_i,$ assign a conjugacy class $c_i$ in $\SU$  
$-2 \le c_i \le 2$
and let
$$
{\mathcal C} = \{c_1, c_2, ..., c_n \}.
$$

\begin{defin}~\label{def:1.1}
The relative character variety with respect to ${\mathcal C}$ is
$$
{\mathcal M}_{\mathcal C} = \{[\rho] \in {\mathcal M} : \tr(\rho(C_i)) = 
c_i, 1 \le i \le n \}.
$$
\end{defin}
The space ${\mathcal M}_{\mathcal C}$ is compact, but possibly singular.
The set of smooth points of ${\mathcal M}_{\mathcal C}$ possesses a natural
symplectic structure which gives rise 
to a finite measure $\mu$ on ${\mathcal M}_{\mathcal C}$ (see \cite{Go1, Go2}).

Let  $\o{Diff}(M, \partial M)$ be the group of diffeomorphisms fixing
$\partial M$.  The mapping class  group $\Gamma$ is defined
to be $\pi_0(\o{Diff}(M, \partial M))$.  The group $\Gamma$
acts on $\pi_1(M)$ fixing the $C_i$'s.
It is known that $\mu$
is invariant with respect to the $\Gamma$-action. In \cite{Go1},
Goldman showed that, with respect to the measure $\mu$,
\begin{thm} [Goldman] \label{thm:1}
The mapping class group $\Gamma$ acts
ergodically on ${\mathcal M}_{\mathcal C}$.
\end{thm}
Goldman also showed that the mapping class action 
is weak-mixing, but not mixing.

Since ${\mathcal M}_{\mathcal C}$ is a variety, one may also study the
topological dynamics of the mapping class group action.  
The topological-dynamical problem is considerably more delicate.  
To begin with,
not all orbits are dense in ${\mathcal M}_{\mathcal C}$.  If 
$\rho \in \Hom(\pi_1(M),G)$
where $G$ is a proper closed subgroup of $\SU$ and $\gamma \in \Gamma$, then
$\gamma(\rho) \in \Hom(\pi_1(M),G)$.   In other words, 
$\Hom(\pi_1(M),G)/G
\subset {\mathcal M}_{\mathcal C}$ is invariant with respect to the 
$\Gamma$-action.
The case of the one-holed torus has been dealt with in \cite{Pr1}:
\begin{thm} \label{thm:2}
Suppose that $M$ is a torus with one boundary component and
$\rho \in \Hom(\pi_1(M),\SU)$ such that $\rho(\pi_1(M))$
is dense in $\SU$.
Then the $\Gamma$-orbit of the conjugacy class
$[\rho] \in {\mathcal M}_{\mathcal C}$
is dense in ${\mathcal M}_{\mathcal C}$.
\end{thm}
This paper deals with the general case of $g>0$ and proves:
\begin{thm} \label{thm:3}
Suppose that $M$ is a Riemann surface with boundary 
having genus greater than zero.  Let
$\rho \in \Hom(\pi_1(M),\SU)$ such that $\rho(\pi_1(M))$
is dense in $\SU$.
Then the $\Gamma$-orbit of the conjugacy class
$[\rho] \in {\mathcal M}_{\mathcal C}$
is dense in ${\mathcal M}_{\mathcal C}$.
\end{thm}


The group $\SU$  double covers the group $\SO$:
$$
p : \SU \longrightarrow \SO.
$$
The group $\o{SO}(3)$ contains subgroups isomorphic to:
$\o{O}(2)$, as well as
the symmetry groups of the regular polyhedra: $T'$ (the tetrahedron),
$C'$ (the cube), $D'$ (the dodecahedron), and their subgroups.
The inverse images of $\o{O}(2)$, $T'$, $C'$, $D'$ by the projection
$p$ are called $\Pin$, $T$, $C$, $D$, respectively.  The identity
component of $\Pin$ is called $\Spin$.
Let  $\rho \in \Hom(\pi_1(M),\SU)$.  
Theorem~\ref{thm:3} implies that if $\rho(\pi_1(M))$ is not contained in 
a group isomorphic to $C, D,$ or $\Pin$, then the $\Gamma$-orbit 
of the conjugacy class $[\rho] \in {\mathcal M}_{\mathcal C}$
is dense in ${\mathcal M}_{\mathcal C}$.

Theorem~\ref{thm:3} covers all moduli spaces except that
of the $n$-holed sphere.  

\begin{conjecture} Theorem \ref{thm:3} holds for all
Riemann surfaces with boundary. 
\end{conjecture}

\noindent To prove this conjecture, 
one must first carry out an analysis of the four-holed sphere
similar to the one performed on the one-holed torus in 
 \cite {Pr1}.  The argument in \cite {Pr1} involves a 
detailed combinatorial study of certain four-term
 trigonometric Diophantine equations (see \cite{Co1}). 
For the four-holed sphere, the analysis would involve 
trigonometric Diophantine equations
with 6 terms. Provided such a result can be obtained, the 
techniques provided in this paper could be adapted to
extend such a four-holed sphere result
to the $n$-holed sphere.

\subsection{Outline of the Proof}
A pants decomposition $\mathcal P$ of $M$ gives rise to a smooth
open dense subset ${\mathcal M}_{\mathcal P}$ that is an integrable system
inside
the moduli space ${\mathcal M}_{\mathcal C}.$ 
Hence, one obtains the following diagram \cite{Go1}:
$$
\begin{CD}
{\mathcal M}_{\mathcal P}        @>f_{\mathcal P}>>       P\\
@VV{i}V                                 @VV{i}V\\
{\mathcal M}_{\mathcal C}      @>f_{\mathcal P}>>    P'
\end{CD}
$$
where ${\mathcal M}_{\mathcal P}$ is a torus bundle over $P$ and
$P' \subset [-2,2]^N$, 
$N = \frac 1 2 \dim({\mathcal M}_{\mathcal C}).$
The subgroup $\Gamma_{\mathcal P} \subset \Gamma$ that preserves the fibres of 
$f_{\mathcal P}$ acts as rotations on each fibre with angles depending
on the base coordinates \cite{Go1}.
Section~\ref{sec:integrable} gives a brief outline of 
this integrable system and the decomposition of the $\Gamma$-action.

The proof of Theorem~\ref{thm:3} involves two steps.
Suppose that
$[\rho]$ is a generic representation (i.e.
$\rho(\pi_1(M))$ is dense in $\SU$).  Let $\Gamma([\rho])$
denote the $\Gamma$-orbit of $[\rho]$.  
The first step is to
show that if $f_{\mathcal P}(\Gamma([\rho]))$ is dense
in $P$, then $\Gamma([\rho])$ is dense in ${\mathcal M}_{\mathcal C}$
(Corollary~\ref{cor:densebase}).
The second step involves proving the base density theorem, i.e., the density of 
$f_{\mathcal P}(\Gamma([\rho]))$ in $P$.

As the problem deals with arbitrary genus, the proof necessarily
involves induction, with the one-holed torus as the base case.  
To get things started, 
for a generic representation $\rho,$
one first shows that there is a one-holed torus $T$ inside $M$ such
that the restriction of $\rho$ to $\pi_1(T)$ is generic.
This is a detailed combinatorial calculation which is outlined in
Section~\ref{sec:handle} and carried out in the
Appendix.

After obtaining a generic handle, we proceed to demonstrate the
base density theorem for the $(n+2g-2)$-holed torus.  
An analysis of the case of the four-holed sphere is 
required to get the induction process started,
and is used to prove the result for the case of the two-holed torus.
From there, the case of the three-holed torus is proven which, in turn, 
is used to prove the case of the $(n+2g-2)$-holed torus.

To complete the proof, the $2g-2$ holes of the $(n+2g-2)$-holed torus
are grouped in pairs and each pair is glued along their boundary
to obtain the original surface $M$ with genus $g$ and $n$ boundary
components.  
Section~\ref{sec:induction} completes the induction process.

\subsection{Some definitions}
Fix a surface $M$ with genus $g > 0$ and $n$ boundary components.
Then $M$ may be described as a $2g$-gon with $n$ holes in the middle,
with appropriate identifications.  More precisely,
the fundamental group $\pi_1(M,O)$ is generated by
$S = \{A_i\}^{2g+n}_{i=1},$
subject to the relation
$$
(\prod^{g}_{i=1} [A_i,A_{i+g}]) (\prod^{2g+n}_{i=2g+1} A_i) = e.
$$
\begin{defin}
\mbox{}

\begin{enumerate}
\item A representation $\rho$ into $\SU$ is generic
if the image of $\rho$ is dense in $\SU$.
\item A handle \cite{Ga1} $(A,B)$ consists of two simple
loops $A,B \in \pi_1(M,O)$
crossing at $O$, but otherwise disjoint.
\item Suppose $G \subset \SU$.  A representation $\rho$ is said
(resp. not) to be $G$ if $Im(\rho)$ is
(resp. not) contained in some (resp. any) isomorphic copy of $G$ in $\SU$.
\item Associated to each simple loop $A \in \pi_1(M,O)$ is the 
Dehn twist in $A$ represented in $\Gamma$ by
a diffeomorphism of $M$ supported on a tubular neighborhood
of $A$ in $M$ (an annulus). The action of the Dehn twist
amounts to cutting $M$
at $A$ and twisting one of the resulting boundary circles once in the 
direction of $A$ and re-identifying the two circles.  
\item With a fixed representation $\rho,$ $X\in \pi_1(M,O),$
 and $\gamma \in \Gamma$,
we write $X$ for $\rho(X)$ and $\gamma(X)$ for $\gamma(\rho)(X)$
when there is no ambiguity.  A small letter will be used to denote the
trace of the matrix represented by the corresponding
capital letter.  For
example, we use $x$ to  denote $\tr(\rho(X))$.  In this setting,
$\gamma(x)$ denotes $\tr(\gamma(\rho)(X))$.
\item  Let $\langle V, {\rm d} \rangle$ be a metric space.
For $\epsilon > 0$, a set $U$ is
$\epsilon$-dense in $V$ if, for each $p \in V,$
there exists a point $q \in U$ such that
$
0 < {\rm d}(p, q) < \epsilon.
$
If $U=V$, we simply say that $V$ is $\epsilon$-dense.
\end{enumerate}
\end{defin}

\centerline{\sc Acknowledgments}
 
Eugene Xia was at the University of Arizona when
much of this research was carried out.
He also thanks IH\'{E}S for its hospitality.
We also mention here that William Goldman has
communicated to us that Michael Kapovich has 
recently obtained similar results.

\section{The Moduli Spaces as Tori Bundles}\label{sec:integrable}
We begin by giving a brief description of the integrable system on the
moduli space \cite{Go2}.  
Suppose $M$ has genus $g \ge 1$ and $n \ge 0$ boundary components.
Let 
$$
\mathcal C = \{c_1,...,c_n \}
$$
be a fixed set of conjugacy classes with 
the first 
$m$ classes not equal to $\pm 2$ and
$c_{n-m+1}, ..., c_n \in \{\pm 2\}.$
Then the real dimension of ${\mathcal M}_{\mathcal C}$ is
$6g-6+2m.$ 
Since the case of the torus is well understood,
we assume throughout the remainder of the paper
that $g>1$ or $m>0.$

\subsection{Pants decompositions}
There is a map
$$
f_{\mathcal P} : {\mathcal M}_{\mathcal C} 
\longrightarrow [-2,2]^{N}
$$
that arises from a so-called pants decomposition ${\mathcal P}$ of $M$,
where $N= 3g-3 +m$.  
For detailed
information, see \cite{Go1}.  
The idea is that, by cutting along
$3g-3+m$ circles on $M$,
the surface $M$ can
be decomposed into $2g-2+m-1$ three-holed spheres
and one exceptional $(n-m +3)$-holed sphere with $n-m$ 
of its boundary components assigned $\pm 2$.

\

\centerline{\epsffile{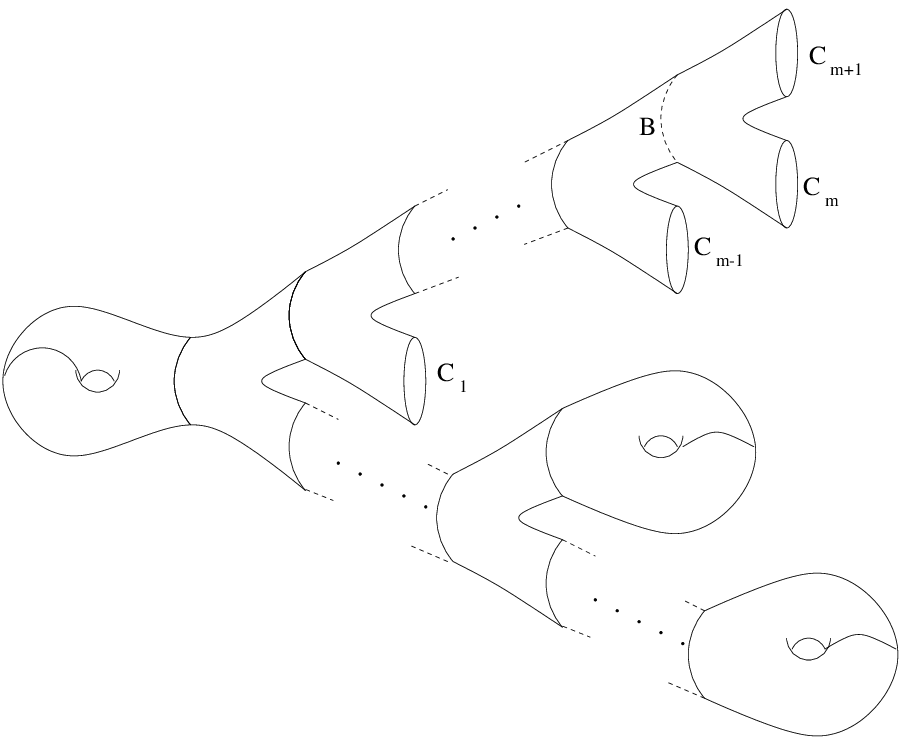}}
\centerline{F{\sc igure} 1: Decomposition of $M$ for $n = m + 1$.}

\

Figure 1 demonstrates such a decomposition when $n = m + 1$, i.e.,
$c_i \neq \pm 2,$ for $1 \le i \le m,$ and  $c_{m+1}=\pm 2$. In this 
specific case, 
we do not cut $M$ along $B,$ thus obtaining 
$2g-2+m-1$ three-holed spheres and one four-holed sphere.
Now in such a case,
if $c_{m+1} = 2$, then the boundary $C_{m+1}$ effectively does not
exist, both with respect to representations and the mapping class
group action.  On the other hand, if $c_{m+1} = -2$ and $n > 1$, then
one can simply consider the moduli space of genus $g$ having $n-1$
boundary components with $\mathcal C = (c_1, ..., -c_m)$.
\begin{rem} \label{rem:reduction} For $ m \ge 1,$
 the moduli space of a genus $g$ surface $M$ with $n$ boundary components
and $\mathcal C = (c_1, ..., c_m, \pm 2, ..., \pm 2)$ 
can be identified with 
the moduli space of a genus $g$ surface $M'$ with $m$ boundary
components and $\mathcal C' = (c_1, ..., \pm c_m).$ 
\end{rem}

In light of Remark~\ref{rem:reduction}, 
we will assume that either
$C_i \neq \pm I$ for all $i$ or that  $n=1$ and
$C_1 = \pm I.$ 

Fix a decomposition $\mathcal P$ of $M$ as described above.
This provides $3g-3+m$ loops
$B_1,...,B_{3g-3+m} \in \pi_1(M,O).$
Let $[\rho] \in {\mathcal M}_{\mathcal C}$ be given.
Then
$$
f_{\mathcal P}([\rho]) = (b_1,...,b_{3g-3+m})
$$
is the desired map, where $b_i=\tr(\rho(B_i))$.
Let $\beta = (b_1, ..., b_{3g-3+m}).$

\subsection{The integrable system $({\mathcal M}_{\mathcal P}, f_{\mathcal P})$}
Let $P'$ be the image of $f_{\mathcal P}$ and let 
$P = P' \setminus \partial P'$.  
The map $f_{\mathcal P}$ restricted to ${\mathcal M}_{\mathcal P}  =
f_{\mathcal P}^{-1}(P)$ is a 
submersion \cite{Go1}.
Denote by $\Gamma_{\mathcal P} \subset \Gamma$ the 
stabilizer of the fibres of $f_{\mathcal P}$.

\begin{prop} \label{prop:integrable}
The set $P'$ consists of all
$\beta \in [-2,2]^{3g-3+m}$
that simultaneously satisfy the $2g-2+m-1$ inequalities 
$$b_i^2+b_j^2+b_k^2 - b_ib_jb_k \le 4,$$
where the three curves 
$B_i, B_j$ and $B_k$ bound a triply punctured sphere in the 
decomposition $\mathcal P$
of $M$ (this includes the possibility of $B_i$ being a
boundary curve in $\partial M$).  In addition, if
$B_i, B_j, B_k$ are the three boundaries of
the exceptional $(n-m+3)$-holed sphere with 
$b_i, b_j, b_k \neq \pm 2$, then 
$$b_i^2+b_j^2+b_k^2 - c b_ib_jb_k \le 4,$$ 
where $c$ is the sign of $\Pi_{i = m+1}^n C_i$.
Suppose $\beta \in P$.
Then there is a $\Gamma_{\mathcal P}-$equivariant homeomorphism
$$h: f_{\mathcal P}^{-1}(\beta) \rightarrow T^{3g-3+m}$$
such that for all $\xi \in f_{\mathcal P}^{-1}(\beta)$ and
$(n_1, ..., n_{3g-3+m})\in \mathbb Z^{3g-3+m}$
$$h: \tau_1^{n_1}...\tau_{3g-3+m}^{n_{3g-3+m}} \xi \mapsto
\left[
\begin{array}{c}
e^{in_1\theta_1} h_1\\
\vdots \\
 e^{in_{3g-3+m}\theta_{3g-3+m}} h_{3g-3+m} \\
\end{array}
\right],
$$ 
where
$$
h(\xi) =
\left[
\begin{array}{c}
h_1\\
\vdots \\
h_{3g-3+m} \\
\end{array}
\right],
$$
$\theta_j = \cos^{-1}(b_j/2),$ and $\tau_i$ is the
action of the Dehn twist in $B_i.$
\end{prop}
\begin{proof}
See \cite{Go1}.
 
\end{proof}

In short, $({\mathcal M}_{\mathcal P}, f_{\mathcal P})$ is
an integrable system.
The real dimensions of $P$ and ${\mathcal M}_{\mathcal P}$ are $3g-3+m$
and $6g-6+2m$, respectively.  
Denote by ${\mathcal M}_{\mathcal C}^s$ the subset of
irreducible representations of ${\mathcal M}_{\mathcal C}$.
\begin{lemma}
${\mathcal M}_{\mathcal C}^s$ is smooth, open, and dense in
${\mathcal M}_{\mathcal C}$.
\end{lemma}
\begin{proof}
The space ${\mathcal M}_{\mathcal C}$ corresponds to the
moduli space ${\mathcal M}_{\mathcal C}^{ss}$ of semi-stable parabolic $\SU$-bundles on $M$.
Moreover
${\mathcal M}_{\mathcal C}^s$ corresponds to the stable
subspace and is open and dense in ${\mathcal M}_{\mathcal C}^{ss}$ \cite{Se1}.
\end{proof}
\begin{prop} \label{prop:dense}
The subset ${\mathcal M}_{\mathcal P}$ is open and dense in
${\mathcal M}_{\mathcal C}^s$.
\end{prop}
\begin{proof}
A direct calculation from Proposition~\ref{prop:integrable}
shows that ${\mathcal M}_{\mathcal C}^s \setminus {\mathcal M}_{\mathcal P}$ is a
real algebraic subvariety with positive co-dimension.
The result then follows from the fact that ${\mathcal M}_{\mathcal C}^s$
is smooth and has dimension $6g - 6 + 2m$.
\end{proof}
Together Proposition~\ref{prop:integrable} and ~\ref{prop:dense} imply:
\begin{cor} \label{cor:densebase}
Let $[\rho] \in {\mathcal M}_{\mathcal C}$ and $\Gamma([\rho])$
be the $\Gamma$-orbit of $[\rho]$.  If 
$f_{\mathcal P}(\Gamma([\rho]))$ is dense in $P$,
then $\Gamma([\rho])$ is dense in ${\mathcal M}_{\mathcal C}$.
\end{cor}

\section{Generic Representations and Handles}\label{sec:handle}

Here we adapt an idea in \cite{Ga1} to find a generic handle (one-holed
torus) inside $M$.  Let $\rho \in \Hom(\pi_1(M,O),\SU).$ 
We first observe that moving the base point $O$, has the effect of
conjugating $\rho$ by an element in $\SU$.

\begin{prop} \label{prop:generic}
For any generic $\rho \in \Hom(\pi_1(M,O),\SU)$, there exists a handle
$(A, B)$ such that $\rho\vert_{\langle A, B \rangle}$ is
generic.
\end{prop}

The proof of Proposition~\ref{prop:generic} is highly computational.
Therefore, we outline its general structure here and leave the details
to the Appendix. One should first consult Section~\ref{sec:1hole} and
\cite{Pr1}, as the proof involves the moduli spaces of the one-holed torus and
uses many ideas pertaining to those spaces.

Suppose $\rho$ is a generic
representation.  For any $i \le g$,
the group $\langle A_i, A_{g+i} \rangle$ may very well be contained
in a proper closed subgroup $G \subset \SU$.  However, the crucial observation
is that for any $1 \le j \le 2g+m$  with $j \neq i$ and $j\neq g+i$, the following
are also handles: 
 $(A_i, A_{g+i} A_j)$,  $(A_i A_j, A_{g+i}),$
 $(A_i,  A_jA_{g+i})$,  $(A_j A_i, A_{g+i})$ (see Figures 2 and 3).
The fact that $\rho$ is generic implies that there exists
$j$ such that $A_j \not\in G$.
 

\

\centerline{\epsffile{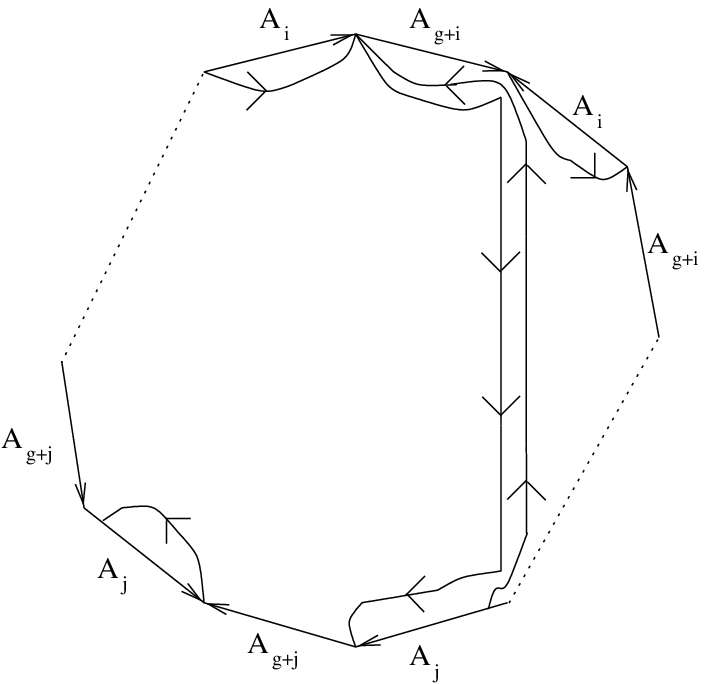}}
\centerline{F{\sc igure} 2: $(A_i, A_{g+i}A_j)$ is a handle.}
 
\ 

\centerline{\epsffile{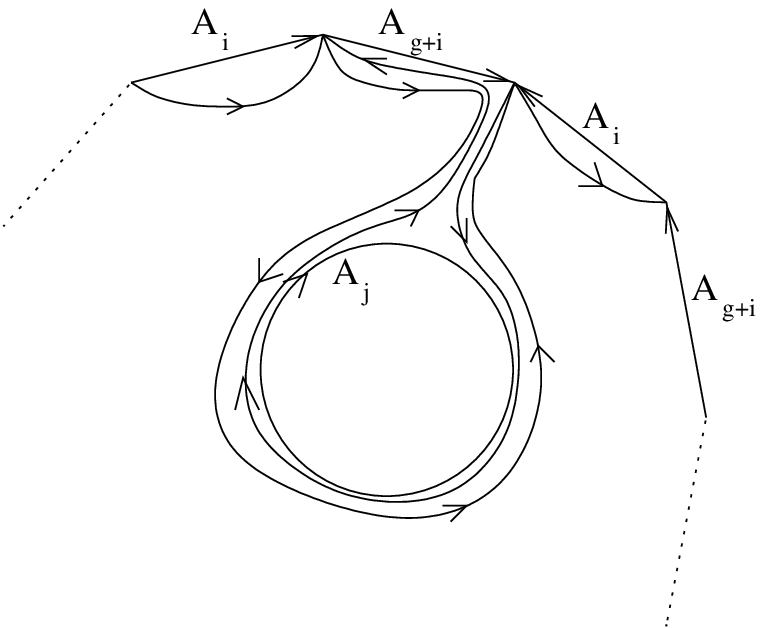}}
\centerline{F{\sc igure} 3: $(A_i,  A_{g+i}A_j)$ is a handle.}
 
\
 
One begins by assuming that $\langle A_i, A_{g+i} \rangle$ is
$\Spin$ and then shows that there exists $j$ such that
$\langle A_i, A_{g+i} A_j \rangle$ is not $\Spin$.
From there, one moves onto the groups $\Pin,$ $T,$ $C,$ and finally to $D$.
Note that the 
computational complexity increases drastically from one case
to the next, not only because the groups are more complicated, 
 but also because it is a recursive process.  For example, suppose
that $\langle A_i, A_{g+i} \rangle $ is $D$ and one
wants to show that there exists $j$ such
that $\langle A_i, A_{g+i} A_j \rangle$
is generic.   One must not only show that $\langle A_i, A_{g+i} A_j \rangle$
is not
$D$, but in addition that it is neither $C$ nor $\Pin$.
 
\section{The Three-holed Sphere}  
Suppose $M$ is a  three-holed sphere.
Then $\pi_1(M)$ has a presentation:
$$
\langle A, B, C: ABC = I \rangle,
$$
where $A, B,$ and  $C$ represent the homotopy classes of the three boundaries
of $M$.  
\begin{prop} \label{prop:Pinon3hole}
\mbox{}

\begin{enumerate}
\item A representation $\rho$ on a three-holed sphere is a 
$\Spin$-representation
if and only if $a^2+b^2+c^2-abc-4=0$.
\item A representation $\rho$ on a three-holed sphere is
$\Pin$ and not
$\Spin$ if and only if 
$a^2+b^2+c^2-abc-4\neq 0$
and 
at least two of the three: 
$A,$ $B,$ $AB,$
have zero trace.
\end{enumerate}
\end{prop}
\begin{proof}
If $\rho$ is a $\Spin$ representation, then up to conjugation:
$$\rho(A) = \left[
\begin{array}{cc}
\cos\theta & \sin\theta \\
-\sin\theta & \cos\theta \\
\end{array}
\right], \ \ 
\rho(B) = \left[
\begin{array}{cc}
\cos\phi & \sin\phi \\
-\sin\phi & \cos\phi \\
\end{array}
\right], \ \ $$
$$\rho(AB) = \left[
\begin{array}{cc}
\cos(\theta+\phi) & \sin(\theta+\phi) \\
-\sin(\theta+\phi) & \cos(\theta+\phi) \\
\end{array}
\right]. \ \ $$

Note that: $$a^2+b^2+c^2-abc-4=$$
$$4\cos^2\theta+4\cos^2\phi+4\cos^2(\theta+\phi)-8\cos\theta\cos\phi\cos(\theta+\phi)-4=0.$$

The other direction follows from the uniqueness of characters \cite{Go1}.
This proves (1).

Suppose $\rho$ is $\Pin$, but not
$\Spin,$ then at least two of the following 
$A,$ $B,$ and $AB$
are in $\Spinm$.
Since $A \in \Spinm$ implies $\tr(A) = 0$, at least two of the
three global coordinates of $[\rho]$ must be zero.  See \cite{Pr1}
for a similar proof in the case of the one-holed torus.
The other direction follows similarly from the uniqueness of 
characters \cite{Go1}.  This proves (2).
\end{proof}

\section{The One-holed Torus} \label{sec:1hole}
We briefly summarize some relevant results that appear in 
\cite{Go1} and \cite{Pr1}.
Suppose that $M$ is a one-holed torus.  
The fundamental group $\pi_1(M)$ has a presentation
$$
\pi_1(M) = \langle X,Y,K | K = XYX^{-1}Y^{-1} \rangle,
$$
where $K$ represents the element corresponding to the boundary component
as in Figure 4.

\

\centerline{\epsffile{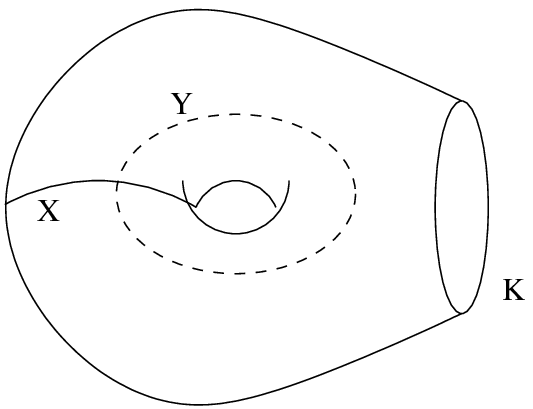}}
\centerline{F{\sc igure} 4: The  one-holed torus}

\


Let
$$
E = \Hom(\pi_1(M), \SU)/\SU.
$$
A representation class $[\rho] \in E$ is determined by
$$
x = \tr(\rho(X)), y = \tr(\rho(Y)), z = \tr(\rho(X Y)).
$$
This provides a global coordinate chart:
$$
[\rho] \stackrel{F}{\longmapsto} (\tr(\rho(X)),
\tr(\rho(Y)), \tr(\rho(X Y))).
$$
In addition, $k = \tr(\rho(K))$ is given by the formula
$$ 
k = \tr(\rho(K)) = x^2 + y^2 + z^2 - x y z - 2.
$$
Let
$$
E_k = \{(x,y,z) \in [-2,2]^3: x^2 + y^2 + z^2 - x y z - 2 = k\},
$$
then
$$
E = \bigcup_{k \in [-2,2]} E_k.
$$
For $ -2 < k < 2$, the set $E_k$ is a smooth two-sphere,  
the set $E_2$ is a singular sphere, and $E_{-2} =(0,0,0)$.

The mapping class group $\Gamma$ is generated by the maps
$\tau_X$ and $\tau_Y$ induced by the Dehn twists
in $X$ and $Y$, respectively.
With respect to the global coordinate, the action can be described
explicitly \cite{Go1}:
$$
\tau_X(x,y,z) = (x, z, xz - y)
$$
$$
\tau_Y(x,y,z) = (z, y, yz - x).
$$
The action of $\tau_X$ fixes $x$ and $k$, and preserves
the ellipse
$$
X_k(x) = \{x\} \times \{ (y,z) : \frac{2+x}{4} (y+z)^2 + \frac{2-x}{4} (y-z)^2
= 2 + k
- x^2 \}.
$$

A change of coordinates transforms $X_k(x)$ into the circle
$$
X_k(x) = \{x\} \times \{(\tilde{y},\tilde{z}) :
\tilde{y}^2 + \tilde{z}^2 = 2 + k - x^2 \}
$$
(see  \cite{Go1, Pr1}).
In this new coordinate system, $\tau_X$ acts as a rotation
by $\cos^{-1}(x/2)$.  In short, the sphere $E_k$ is
the union of circles
$$
E_k = \bigcup_{x} X_k(x) ,
$$
and $\tau_X$ rotates (up to a coordinate transformation)
each level set $X_k(x)$ by an angle of
$\cos^{-1}(x/2)$.
Similarly, there is a coordinate transformation
so that 
$\tau_Y$ acts as a rotation of $Y_k(y)$ by an angle of
$\cos^{-1}(y/2).$

\begin{prop} \label{prop:pin}
The space of $\Spin$ representation classes
consists precisely of $E_2$. The
$\Pin$ representation classes consist of $E_2$ and the intersections of
the three coordinate axes with $E$.
For each $-2 < k < 2$, there are exactly
six points corresponding to $\Pin$ representation classes
in $E_k$. Moreover,
a representation class  $(x, y, z)\in E_k$ with $-2 < k < 2$
and $x \neq 0$
is $\Pin$ if and only if  $k=x^2-2.$
\end{prop}
\begin{proof}
See \cite{Pr1}.
\end{proof}

\begin{rem} \label {rem:rotation1}
This is the first explicit example of Proposition~\ref{prop:integrable} with
$g=1$ and $m=1$.  One obtains a pair of pants by cutting along $X$ (resp. $Y$).
The important property is that if $X_k(x)$ (resp. $Y_k(y)$) is a non-degenerate
circle, 
then $\tau_X$ (resp. $\tau_Y$) acts on the fibre $X_k(x)$ (resp. $Y_k(y)$)
as a rotation with an angle depending solely on $x$ (resp. $y$) and 
independent of either $k$ or $y$ (resp. $x$).  
In particular, if a representation $\rho$ is not 
$\Pin$, then neither $\tau_X$ nor
$\tau_Y$ fix $[\rho]$.
\end{rem}


\begin{rem} \label{rem:generic}
Let $\rho$ be generic. By Theorem~\ref{thm:2}, the 
$\langle \tau_X, \tau_Y \rangle$-orbit
$\O$ of $[\rho]$ is dense in $E_k$.  Hence, there is a number 
$r > 0$ such that the set $\{(x,y) : (x,y,z) \in \O\}$ is dense
in $R = [-r,r]^2$.  
In particular, by Dehn twisting in $\tau_X$ and $\tau_Y$,
one can always assume that $x$ and $y$ are simultaneously
off of any finite set of values
and are arbitrarily close to zero.  
\end{rem}

\begin{rem} \label{rem:1hole} 
In addition to Theorem~\ref{thm:2}, it was shown in \cite{Pr1} that
the $k$-values for the surjective $T, C, D$ representations are $0,
\frac{1 \pm \sqrt{5}}{2}, 1$.
We define the set of special values:
$$
{\mathcal S} = \{\frac{1 \pm \sqrt{5}}{2}, 0, 1 \}
$$
In particular, if $(x,y,z)$ is not $\Pin$ and  $k \not\in {\mathcal S}$,
then $(x,y,z)$ is generic.
\end{rem}

\section{The Four-Holed Sphere}

We first review some results that appear
in \cite{Be1} and \cite{Go1}.
Suppose $M$ is a four-holed sphere.  Then the fundamental
group $\pi_1(M,O)$ admits a presentation
$$
\langle A, B, C, D : ABCD = I \rangle.
$$
\subsection{The moduli space $\M$}
The moduli space $\M$ for the four-holed sphere is six-dimensional.  
Let $g : \M \longrightarrow [-2,2]^4$ be the map 
defined by $g([\rho]) = (a, b, c, d)$.
\

\centerline{\epsffile{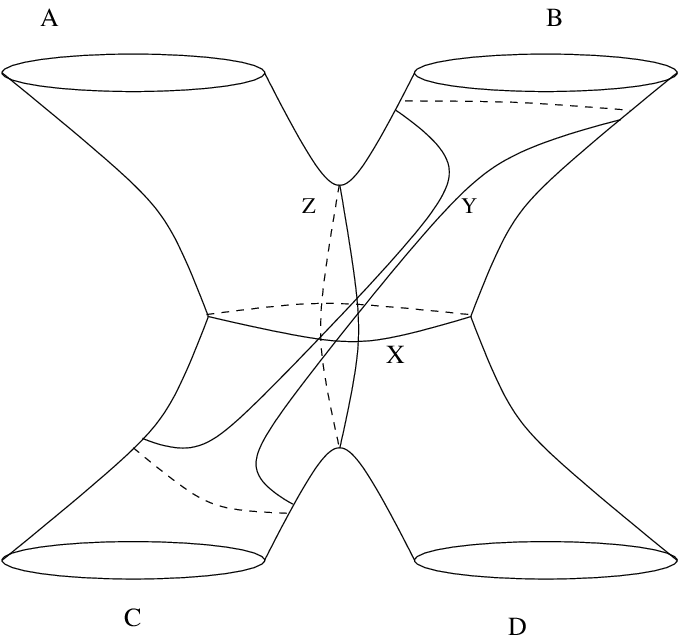}}
\centerline{F{\sc igure} 5: The  four-holed sphere}

\

Set $X = AB, Y = BC,$ and  $Z = CA$ (see Figure 5).
Then the set $E_{\kappa} = g^{-1}(a,b,c,d)$
satisfies the equation
$$
x^2 + y^2 + z^2 + xyz= (ab+cd)x + (ad + bc)y + (ac+bd)z 
-(a^2 + b^2 + c^2 + d^2 + abcd -4).
$$

Let 
$$
I_{a,b}= \bigg [\frac { ab - \sqrt{(a^2-4)(b^2-4)}}2, \frac {ab + \sqrt{(a^2-4)
(b^2-4)}}2 \bigg ].
$$
For any $\kappa =(a,b,c,d) \in [-2, 2]^4 $ with
$I_{a,b} \cap I_{c,d} \neq \emptyset$,
the $x$-level sets $X_{\kappa}(x) \subset E_\kappa$  
are ellipses (possibly degenerate)
in $y$ and $z$ given by:
$$
\frac{2+x}{4}((y+z) - \frac{(a+b)(d+c)}{2+x})^2 + 
\frac{2-x}{4}((y-z) - \frac{(a-b)(d-c)}{2+x})^2 = 
$$
$$
\frac{(x^2-abx+a^2+b^2-4)(x^2-cdx+c^2+d^2-4)}{4-x^2}.
$$
There are similar descriptions for the $y$- and $z$-level sets  
$Y_\kappa(y)$ and $Z_\kappa(z)$ respectively (see \cite{Go1}). 

For $x$ in the interior of  $I_{a,b} \cap I_{c,d}$
the level set $X_{\kappa}(x)$ is an ellipse centered
at 

\begin{equation} \label{eq:center1}
\left\{
\begin{array}{ll}
y_c(x) & = \frac{2(2(ad+bc)-x(ac+bd))}{4-x^2} \\  
z_c(x) & = \frac{2(2(ac+bd)-x(ad+bc))}{4-x^2}.
\end{array}
\right.
\end{equation}

Figure 6 below displays several 
$Z_\kappa(z) \subset E_\kappa$ for a particular $\kappa.$

\

\centerline{\epsffile{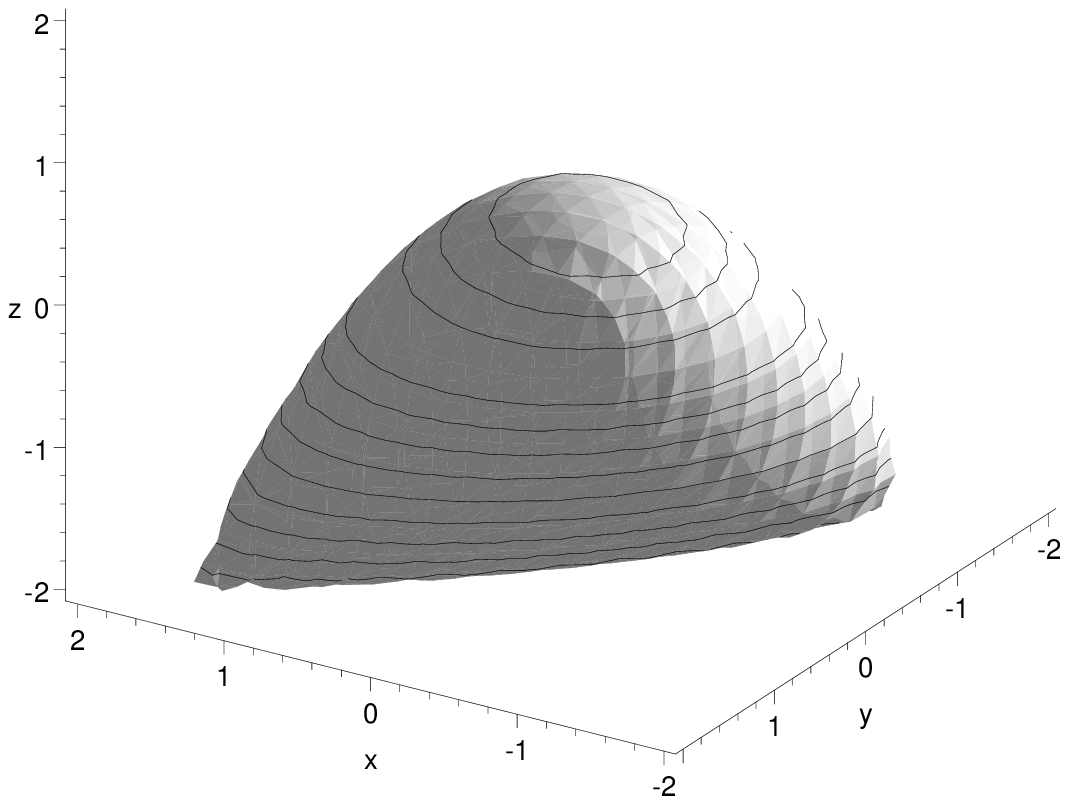}}
\centerline{F{\sc igure} 6: The sphere $E_\kappa$ for $\kappa =(.5, -1, -.2, 1.2).$}

\

For fixed $\kappa = (a,b,c,d)$, the
$x$-coordinates in $E_\kappa$ take on all the values inside 
$I_{a,b} \cap I_{c,d}$ 
(see  Proposition~\ref{prop:integrable} and \cite{Be1}).
In particular,
$$E_\kappa = \bigcup_{x\in I_{a,b} \cap I_{c,d}}X_\kappa(x).$$
By symmetry, similar constructs can be made  
for the $y$- and $z$-coordinates.

\subsection{The mapping class action}
In local coordinates, 
the actions of $\tau_X, \tau_Y, \tau_Z$ are
$$\left[
\begin{array}{c}
y\\
z\\
\end{array}
\right]
\stackrel{\tau_X}{\longmapsto}
\left[
\begin{array}{c}
ad+bc-x(ac+bd-xy-z)-y\\
ac+bd- xy-z\\
\end{array}
\right],
$$
$$\left[
\begin{array}{c}
z\\
x\\
\end{array}
\right]
\stackrel{\tau_Y}{\longmapsto}
\left[
\begin{array}{c}
bd+ca-y(ba+cd-yz-x)-z\\
ba+cd-yz-x\\
\end{array}
\right],
$$
$$\left[
\begin{array}{c}
x\\
y\\
\end{array}
\right]
\stackrel{\tau_Z}{\longmapsto}
\left[
\begin{array}{c}
cd+ab- z(cb+ad-zx-y)-x\\
cb+ad- zx - y\\
\end{array}
\right].
$$
These actions preserve the ellipses $X_\kappa(x) \subset E_\kappa,
Y_\kappa(y) \subset E_\kappa,$ and $Z_\kappa(z) \subset E_\kappa,$
respectively.
After coordinate transformations,
these are rotations by angles $2\cos^{-1}(x/2), 2\cos^{-1}(y/2),$ and
$2\cos^{-1}(z/2)$, respectively \cite{Go1}.
\begin{rem}
This is the second explicit example of Proposition~\ref{prop:integrable} with
$g=0$ and $m=4$.  One obtains two pairs of pants by cutting along $X$
(resp. $Y$ or $Z$).
The important property is that 
if $X_\kappa(x)$ (resp. $Y_\kappa(y)$ or $Z_\kappa(z)$) is a 
non-degenerate circle, then
$\tau_X$ (resp. $\tau_Y$ or $\tau_Z$) acts on the 
fibre $X_\kappa(x)$ (resp. $Y_\kappa(x)$ or $Z_\kappa(z)$)
as a rotation with an angle depending only on $x$ (resp. $y$ or $z$) and 
independent of $\kappa$.
\end{rem}

Let ${\rm d}$ be the metric
on $\R^3$ that generates the box topology.
The metric ${\rm d}$ 
generates the usual topology on $\M_{\mathcal C}$.
For a fixed $\kappa$, the coordinates 
provide an embedding of $E_\kappa$ into $\R^3$ and
$E_\kappa$ inherits the metric $d$ which is described as:
$$
{\rm d}((x_1, y_1, z_1), (x_2, y_2, z_2)) = 
max\{|x_1 - x_2|, |y_1 - y_2|,|z_1 - z_2|\}.
$$

\subsection{Filtration on the level sets} 
We introduce a filtration that is analogous to the one
introduced in \cite{Pr1} for the one-holed torus. 
The Dehn twist $\tau_Y$ acts on the (transformed)
subsets $Y_\kappa(y)$ via a rotation of angle  $2\cos^{-1}(y/2)$.
Thus there is a filtration of the 
$y$-coordinates that yield finite orbits under $\tau_Y$ as
follows: Let $Y_n \subset (-2, 2)$ such that $y \in Y_n$ if and only if
the $\tau_Y$-action on non-fixed points $(x,y,z) \in E_\kappa$ is 
periodic with period less than or equal to $n$.
This gives a filtration
$$
\{0\} = Y_2 \subset Y_3 \subset ... \subset Y_n \subset ...
$$
For example: $Y_2 = \{0\}$, $Y_3 = \{0, 1, -1\}$, 
$Y_4 = \{0, 1, -1, \sqrt 2, -\sqrt 2\}$,
etc. 
Note that the filtration is 
independent of the choice of $\kappa$ and that $Y_n$ is a 
finite set for every $n$. 
By symmetry, there exist similar filtrations $X_n$ and $Z_n,$ with 
$X_n=Y_n=Z_n$ 
as sets.
The following lemmas, though proven for the $Y_n$ filtration, apply 
equally to the other filtrations.

\begin{lemma}\label{lem:path1}
For $\epsilon > 0$ there exists $N(\epsilon) > 0$ 
so that if $y \not \in Y_{N(\epsilon)}$,
then the $\tau_Y$-orbit of $(x,y,z)$ 
is $\epsilon$-dense in $Y_\kappa(y)$ for any $(x,y,z)$ in any $E_\kappa.$
\end{lemma}
\begin{proof}

Since the ellipses $Y_\kappa(y)$ are 
 (possibly degenerate) of uniformly bounded circumferences,
there exists $N(\epsilon) > 1$ such that 
for any $y \not\in Y_{N(\epsilon)}$, the $\tau_Y$-orbit is $\epsilon$-dense
in $Y_\kappa(y)$.  
\end{proof}

Throughout the remainder of the paper,
the moduli spaces of four-holed spheres with
$\kappa =(a,b,c,d),$ having small $|c|$ and $|d|,$ 
will play an important role. We first analyze 
the case of $E_\kappa$ for
$\kappa =(a,b,0,0).$

\begin{lemma} \label{lem:3circle}
Suppose $(a,b,0,0) = \kappa \in [-2,2]^4$.  Then
$X_{\kappa}(x) \subset E_\kappa$ is 
an ellipse (possibly degenerate) centered at $(x,0,0)$.
Thus, $Y_{\kappa}(0)$ 
and $Z_{\kappa}(0)$ intersect every ellipse $X_{\kappa}(x).$ 
\end{lemma}
 
\begin{proof}
This result follows directly from 
 equation (\ref{eq:center1}) which gives
$(y_c(x), z_c(x))$,
the center of the ellipse $X_{\kappa}(x).$ 
\end{proof}
 
\begin{lemma} \label{lem:inband}
Let $a, b \in (-2, 2)$ and $\epsilon >0$.  
Then there exists $\delta > 0$ so that
for any
$\kappa=(a,b,c,d)$ with
$|c|, |d| < \delta,$ 
every ellipse 
$Y_{\kappa}(y) \subset E_{\kappa}$ with
$|y|\le \delta$ has points with $x$-coordinates 
that come
within $\epsilon$
of all possible $x$-coordinates in 
$E_{\kappa}.$
That is, $$\{x: (x,y,z) \in E_{\kappa} \} \subset
\{x \pm \epsilon : (x,y,z) \in Y_{\kappa}(y), \ |y| \le \delta\}.$$

Moreover, each
$X_{\kappa}(x) \subset E_{\kappa}$ satisfies at least one of:
\begin{enumerate}
\item $X_{\kappa}(x)$ has points which realize all $y$- and $z$-coordinates in
$[-\frac {\delta} 2, \frac {\delta} 2].$ 

\item  All $y$- and $z$-coordinates of $X_{\kappa}(x)$ are
inside $[-\delta,\delta].$ 
\end{enumerate} 

\end{lemma}

\begin{proof}
The result holds by the continuous dependence of $E_{\kappa}$
on $c$ and $ d,$ by  the continuous dependence of $X_{\kappa}(x)$
on $x$,
and by the geometry of $E_{\kappa}$
for $\kappa=(a,b,0,0)$ as described in Lemma~\ref{lem:3circle}.
\end{proof}

\section{The Two-Holed Torus}
In this section, we prove Theorem~\ref{thm:3} for the case of
$g=1$ and $m=2$.
The fundamental group $\pi_1(M,O)$ has a presentation
$$
\langle X, Y, A, B | XYX^{-1}Y^{-1}= A B \rangle,
$$
where $A$ and $B$ represent the boundary components.
Let 
$$
K = XYX^{-1}Y^{-1}, W = AX, W' = XB, Z = XY,
$$
(see Figure 7).


\

\centerline{\epsffile{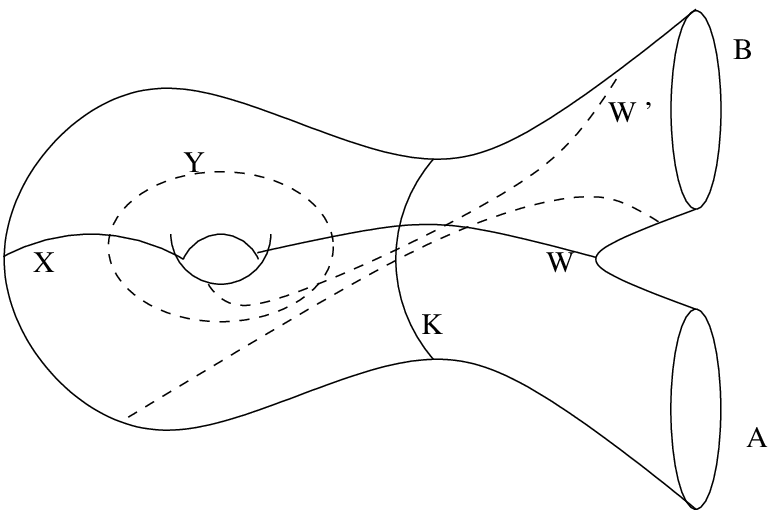}}
\centerline{F{\sc igure} 7: Two-holed torus.}

\

Assume that $\rho \in \Hom(\pi_1(M,O),\SU)$ is
generic.
By Proposition~\ref{prop:generic}, we have that $\langle X, Y \rangle$
is generic.
If we cut $M$ along $X$ (resp. $Y$), then we
obtain a four-holed sphere
with $\kappa=(a,b,x,x)$ (resp., $\kappa=(a,b,y,y)$).
By Remark \ref{rem:reduction} and Theorem \ref{thm:2}, we 
assume that $a, b \neq \pm 2.$ 
For $\kappa = (a,b,x,x)$, the formula for the defining equation
of $E_{\kappa}$ simplifies to:
$$w^2+(w')^2 +k^2 + kww'=$$
$$k(ab+x^2)+wx(a+b)+w'x(a+b)-a^2-b^2-2x^2-abx^2+4,$$ 
and the formulas for the Dehn twists $\tau_K, \tau_W, \tau_{W'}$
become:
$$\left[
\begin{array}{c}
w\\
w'\\
\end{array}
\right]
\stackrel{\tau_K}{\longmapsto}
\left[
\begin{array}{c}
x(a+b) - k(x(a+b) -kw -w')-w\\
x(a+b) -kw -w'\\
\end{array}
\right],
$$
$$\left[
\begin{array}{c}
w'\\
k\\
\end{array}
\right]
\stackrel{\tau_W}{\longmapsto}
\left[
\begin{array}{c}
x(a+b)- w(x^2+ab -ww'-k)-w'\\
x^2+ab -ww'-k\\
\end{array}
\right],
$$
$$\left[
\begin{array}{c}
k\\
w\\
\end{array}
\right]
\stackrel{\tau_{W'}}{\longmapsto}
\left[
\begin{array}{c}
ab+x^2- w'(x(a+b) -w'k-w)-k\\
x(a+b) -w'k-w\\
\end{array}
\right].
$$
Note that $I_{x,x} = [x^2-2, 2].$

\begin{lemma}\label{lem:prelim1}
For the representation $\rho$,
there exists  a sequence of Dehn twists $\gamma \in \Gamma$
so that at least one of the following has non-zero trace:
$\gamma(AX)$, $\gamma(XB)$, $\gamma(AXY),$ $\gamma(XYB),$
$\gamma(AY),$ $\gamma(YB),$ $\gamma(AYX),$  $\gamma(YXB)$, with $\langle \gamma(X),\gamma(Y) \rangle$ generic.
\end{lemma}
\begin{proof}
Suppose that $AX$, $XB$, $AXY,$ $XYB,$
$AY,$ $YB,$ $AYX,$ and  $YXB$ all have zero trace.
If $\tau_K$ preserves $w=w'=0$, then
$x(a+b)=0$ (if not, note that the twist in $\tau_K$
preserves $x, y$ and $k,$ so $\langle \tau_K (X), \tau_K(Y) \rangle$ 
remains generic).
By Remark~\ref{rem:generic}, $x$ can be assumed non-zero.  Hence
$a=-b.$  Then the equation for the four-holed sphere $E_{(a,-a,x,x)}$ 
implies that: 

\begin{equation} \label{eq:0}
k^2+ka^2+2a^2-4=x^2(k-2+a^2).
\end{equation}

By Remark~\ref{rem:generic} and the fact that the left-hand side of 
equation (\ref{eq:0}) is invariant under 
$\langle \tau_X, \tau_Y \rangle,$ it must be that 
$k-2-a^2=0.$
By the formula $\tr(AB^{-1})+\tr(AB)=\tr(A)\tr(B)$ (see \cite{Ma1,Ma2})
and the fact that  $K=AB,$ 
we have that
$\tr(AB^{-1})+2-a^2=-a^2$; this implies that $A=-B$
as matrices.

Suppose $A=-B$. Conjugating by an element in $\SU$, 
we may assume that 
$$X = \left[
\begin{array}{cc}
x_1 + y_1i  & 0\\
0 & x_1 + y_1i\\
\end{array}
\right], \ \ 
Y = \left[
\begin{array}{cc}
x_2 + y_2i  & z_2 \\
-z_2 & x_2 - y_2i\\
\end{array}
\right],\ \ 
$$
and
$$A = \left[
\begin{array}{cc}
a_1 + b_1i  & c_1 + d_1i\\
-c_1 + d_1i & a_1 - b_1i\\
\end{array}
\right]. \ \ 
$$ 
The equations: $\tr(AX) =0,$ $\tr(AY) =0,$ 
$\tr(AXY) =0,$ $\tr(AYX) =0,$ 
together with the matrix equation
$K=-A^{-2}=XYX^{-1}Y^{-1}$
and $\det(A)=\det(X)=\det(Y)=1$ have
solutions that must  include one of the following:
$x_1=0,\pm 1,$ $a_1=0$ 
or $a_1= \pm \sqrt{1-x_1^2}$.  
By Remark~\ref{rem:generic}, one may assume that $a_1 \neq \sqrt{1-x_1^2}$
and $x_1 \neq 0, \pm 1$.
This leaves us with the case $a_1 =0$, which implies that $K=I,$
contradicting the fact that $\langle X, Y \rangle$ is generic.
\end{proof}
Combining Proposition~\ref{prop:pin} and Lemma~\ref{lem:prelim1}, one obtains:
\begin{cor} \label{cor:nonzero}
For a generic $\rho$, we may assume that one of:
$\langle AX, Y \rangle$, 
$\langle XB, Y \rangle$,
$\langle AY, X \rangle$, 
$\langle YB, X \rangle$, 
is not $\Pin.$ 
\end{cor}
By Remark~\ref{rem:rotation1} and
Corollary~\ref{cor:nonzero},  
we assume without loss of generality
that the action $\tau_Y$ does not fix $w=\tr(AX)$.  

\begin{prop} \label{prop:generic2hole}
Suppose $g=1$ and $m=2$ and $\rho$ is a generic representation.
Then the $\Gamma$-orbit $\Gamma([\rho])$ is dense in $\M_{\mathcal C}$.
\end{prop}

\begin{proof} 
Let $\epsilon' > 0$.
Let  $\rho$ be a generic representation. By
Proposition \ref{prop:generic}, $M$ has a 
generic handle $\langle X, Y \rangle$
(we adopt the notation presented in Figure 7).
By Proposition  \ref{prop:dense}, without loss of generality
we may take $[\rho_0] \in \mathcal M_{\mathcal P}$ and show that
there exists $\gamma \in \Gamma$ such that
${\rm d} (\gamma([\rho]), [\rho_0])< \epsilon'.$
Let $[\rho_0] \in \mathcal M_{\mathcal P}$ be a representation class with
$x_0 = \tr(\rho_0(X))$ and $k_0=\tr(\rho_0(K)).$ 
Let $w_0= \tr(\rho_0(W)),$ $w_0'= \tr(\rho_0(W')),$
etc.  

The map $f_{\mathcal P}$ is a submersion on $\M_{\mathcal P}$.
By the implicit function theorem, 
there exists  $\epsilon'>\epsilon>0,$ such that the
$\epsilon$-tubular neighborhood of $f_{\mathcal P}^{-1}(x_0,k_0) \subset 
\mathcal M_{\mathcal C}$ 
looks like
$[x_0-\epsilon, x_0+\epsilon] \times [k_0-\epsilon, k_0+\epsilon]\times T^2.$


Cutting along $K$ and $X$ gives a pants decomposition of $M$.  Hence,
by Corollary~\ref{cor:densebase}, we only need to show that
there exists $\gamma \in \Gamma$ such that $\gamma([\rho])$
satisfies 
$|\gamma(k) - k_0|< \epsilon$ and  $|\gamma(x) - x_0| < \epsilon$.
In other words, the goal is to move the $x$- and $k$-coordinates near
$x_0$ and $k_0$.  The strategy is to use $\tau_Y$ to
first move the $x$-coordinate
near zero so that the moduli space of the four-holed sphere 
(obtained
by cutting along
$X$) contains $k$-coordinates near $k_0$.  At the same time
one needs to ensure that $\tau_W$ has enough points on its orbit
so that the $k$-coordinate can actually be moved
near $k_0$ (this may require a sequence of twists in $\tau_W$
followed by a sequence of twists in $\tau_K$ followed again
by a sequence of twists in $\tau_W$).  This can be accomplished by using the 
generic $\langle X, Y \rangle$ to
move the $y$-coordinate to a sufficiently general position (i.e., outside
of $Y_N$ for some very large $N$).  By Corollary~\ref{cor:nonzero}, this
will produce a $\tau_Y$-orbit in $\langle Y, W \rangle$ with a
sufficiently large number of points.  
Under such conditions, this $\tau_Y$-orbit
will potentially contain points that have a large
number of points on their
$\tau_W$-orbit 
in $E_{(x,x,a,b)}$.
The exceptional case being 
when the $\tau_W$-action fixes $(k,w,w') \in E_{(x,x,a,b)}$.
In this case,
$W_\kappa(w)$ is a point, which implies that
$$
2k = ab+x^2-ww'.
$$
If this occurs, one may try to use the $\tau_{W'}$-action.
However in order to use the $\tau_{W'}$-action, one must
first ensure that $\langle Y, W' \rangle$ itself is not $\Pin.$
Hence under the assumption that $\tau_W$ fixes
$(k,w,w') \in E_{(a,b,x,x)},$ we have the following 
two exceptional cases:

Case 1:  $\langle W', Y\rangle$ is $\Pin$.
This leads to either $y=0$ or
$w'=0$. This implies that either $y=0$ or
\begin{equation} \label{eq:e1}
2k = ab + x^2.
\end{equation}

Case 2: $\langle W', Y\rangle$ is not $\Pin$
but $(k,w,w') \subset E_{(a,b,x,x)}$ 
is fixed by $\tau_{W'}.$ 
Then one must have
$$
\left\{
\begin{array}{lll}
2k & = & ab+x^2-ww' \\[2ex]
2w'& = & x(a+b) - w k \\[2ex]
2w & = & x(a+b) -w' k
\end{array}
\right.
$$
This leads to $2(w'-w)=  k(w'-w).$
Since $k \neq \pm 2$,  we have $w=w'=\frac {x(a+b)}{k+2}.$ 
Moreover, for non-zero $w$, we have $a \neq -b$.
This implies that

\begin{equation} 
x^2 \bigl (1- \bigl(\frac {a+b}{k+2} \bigr )^2 \bigr ) = 2k-ab.
\end{equation}
If $a+b= \pm (k+2)$ and $k = \frac {ab}2,$ 
then  either $a=\pm 2$ or $b= \pm 2.$  This contradicts our assumptions
on $a$ and $b$.
So,
\begin{equation} \label{eq:e2}
x^2 = \frac {2k-ab}{1- \bigl(\frac {a+b}{k+2} \bigr )^2 }.
\end{equation}

Now we choose a generic representation
$\langle X, Y \rangle$
that allows us to avoid the special cases presented above.
The following lemma follows immediately from 
Remark~\ref{rem:generic}, Corollary~\ref{cor:nonzero} and 
the fact that the $\tau_X, \tau_Y$ actions fix $k$.
\begin{lemma} \label{lem:laundry1}
For any integer $J>0$,
there is
$\gamma \in \Gamma$ such that the $\tau_Y$-orbit of $\gamma (\rho)$ has
at least $J$ points 
satisfying the
following conditions:
\begin{enumerate}
\item The $x$-coordinates of these $J$ points have
$|x|$ small enough so that 
$k_0 \pm \frac \epsilon 2 \in I_{a,b} \cap I_{x,x}$ and
$|x| \le \delta,$ where 
$\delta$ is provided in Lemma \ref{lem:inband}
for $E_{(a,b,x,x)}.$ 
\item These $J$ points do not belong to the subvarieties defined by
equations (\ref{eq:e1}), (\ref{eq:e2}),  $y=0$, and $w=0$.
\end{enumerate}
\end{lemma}
Note that Condition 2 of Lemma~\ref{lem:laundry1} ensures that if $\tau_W$ 
fixes $(k,w,w'),$ then $\langle Y, W' \rangle$ is not $\Pin$ and
$\tau_{W'}$ does not fix $(k,w,w')$.

Now choose $\gamma \in \Gamma$ with $J$
sufficiently large in Lemma~\ref{lem:laundry1} so that
one of the $J$ points on the $\tau_Y$-orbit of 
$\gamma(\rho),$ denoted $\gamma_1(\rho),$ 
has a $\tau_W$-orbit (or $\tau_{W'}$-orbit) 
with at least 
$$N(\frac {\delta}{N(\frac  \epsilon 2)+3})+3$$ points with distinct 
$k$-coordinates inside $E_{(a,b,\gamma_1(x),\gamma_1(x))}.$

Since $|\gamma_1(x)| \le \delta$ and $k_0 \pm \epsilon 
\in I_{a,b} \cap I_{\gamma_1(x),\gamma_1(x)},$ 
the sphere  $E_\kappa$ (with $\kappa = (a,b,\gamma_1(x),\gamma_1(x))$),
has points with
$k$-coordinates inside 
$[k_0-\frac \epsilon 2, k_0+\frac \epsilon 2].$ 
Since the
$\tau_W$-orbit (resp. $\tau_{W'}$) of $\gamma_1(\rho)$
has  at least 
$$N(\frac {\delta}{N(\frac  \epsilon 2)+3})+3$$ points with distinct $k$-coordinates, 
one such point
 $\gamma_2(\rho)$ has $k$-coordinate not in $K_{N(\frac {\delta}{N( \frac \epsilon 2)+3})}$ with non-degenerate ellipse
$K_{\kappa}(\gamma_2(k)) \subset E_{\kappa}.$ Thus,
the $\tau_K$ orbit of $\gamma_2(w)$ is 
$\frac {\delta}{N( \frac \epsilon 2)+3}$-dense.
Thus, the $\tau_K$-orbit of  $\gamma_2(\rho)$ has at least 
$N( \frac \epsilon 2)+3$ points with $w$-coordinates 
inside $(-\delta,\delta).$
One such point, $\gamma_3(\rho)$  has $w$-coordinate not
in $W_{N( \frac \epsilon 2)}$ with non-degenerate ellipse
$W_{\kappa}(\gamma_3(w)).$
Thus, the 
 $\tau_W$-orbit of $\gamma_3(\rho)$ is $\frac \epsilon 2$-dense
in $W_\kappa(\gamma_3(w))$. 
This fact, together with the properties of $\delta$ provided in
Lemma \ref{lem:inband}
imply that at least one point, $\gamma_4([\rho])$, in the
$\tau_W$-orbit of $\gamma_3(\rho)$
has $k$-coordinate $\gamma_4(k)$
that comes within $\epsilon $ of $k_0.$

%
%

The one-holed torus $\langle \gamma_4(X), \gamma_4(Y) \rangle$ is
generic 
so long as 
$\gamma_4(k) \not\in \mathcal S$ and $\gamma_4(k) \neq (\gamma_4(x))^2-2$
(see Proposition~\ref{prop:pin} and Remark~\ref{rem:1hole}).
This can be accomplished by replacing $\epsilon$ by $\frac \epsilon {20}$
at the start of the argument.
Since  $\langle  \gamma_4(X),  \gamma_4(Y)\rangle$  is generic, 
we may apply $\tau_X$ and $\tau_Y$
to obtain $\gamma_5\in \Gamma$ so that the $x$-coordinate
$\gamma_5(x)$ of $\gamma_5([\rho])$
is within $\epsilon$ of $x_0.$
Note that both $\tau_X$ and $\tau_Y$
fix $\gamma_4(k)$, thus $\langle \gamma_5(X), \gamma_5(Y)\rangle$
remains generic.
The result now follows by Corollary~\ref{cor:densebase}.

To summarize, one first uses $\tau_Y$ to
obtain points with $x$-coordinates 
near zero having  $\tau_W$ (resp. $\tau_{W'}$) actions 
that generate points having $\tau_K$-actions with a sufficiently
large number of points. 
Next, one uses $\tau_K$-action  
to get the
$w$-coordinate to be near zero on the moduli space of the four holed sphere
obtained by cutting at $X.$
This ensures that one can  move the $k$-coordinate near
$k_0$ (see Lemma \ref{lem:inband}).
%
After these twists, the handle $\langle 
\gamma(X),\gamma(Y) \rangle$ is shown to
remain generic.
Finally move the $x$-coordinate near $x_0$ by using
Theorem \ref{thm:2}.
\end{proof}

The moduli space for the two-holed torus is 4-dimensional.  
By Proposition~\ref{prop:generic2hole}, we obtain the following (analogous
to Remark~\ref{rem:generic}):
\begin{rem} \label{rem:generic2}
Suppose $\rho \in \M_{\mathcal C}$ is generic with coordinates $(x,y,k,w)$.
Then the orbit $\Gamma([\rho])$ is dense in 
$\M_{\mathcal C}$.  
Hence there is a number
$r > 0$ such that the set 
$\{(x',y',k',w') : (x',y',k',w') \in \Gamma([\rho])\}$ is dense 
in a ball of radius $r$ centered at $(x,y,k,w)$.
In particular, by Dehn twisting in $\tau_X, \tau_Y, \tau_K, \tau_W$,
one can always assume that $(x, y, k, w)$ does not belong to 
any proper subvariety of $\M_{\mathcal C}$.
\end{rem}

\section{The $m$-holed torus}

As $m>2$ throughout this section, we assume 
that $C_i \neq \pm I$ for all boundary curves $C_i \in \partial M.$
\begin{prop} \label{prop:noI2}
Let  $M$ be an $m-$holed torus and $\rho$ a generic representation with
generic handle $\langle X, Y \rangle,$ which is part of a
pants decomposition $\mathcal P,$ i.e.,
$K=XYX^{-1}Y^{-1}, X \in \mathcal P .$ Then
there is $\gamma \in \Gamma$
 so that
$\gamma(P) \neq \pm I$ for all $P \in \mathcal P.$
Moreover, the handle $\langle \gamma(X), \gamma(Y) \rangle$ is generic.
\end{prop}

\begin{proof}
 
We first treat the case of $m=3$.
Let $C$ and $D$ be two boundary loops separated from
 $\langle X, Y \rangle$
by
$B=CD.$ Let $A$ be the remaining boundary loop 
(see Figure 8).  Since $\langle X, Y \rangle$ is generic,
we have that $K, X \neq \pm I.$ 
The goal, therefore,  is to find $\gamma \in \Gamma$
such that $\gamma(B) \neq \pm I$.

\
 
\centerline{\epsffile{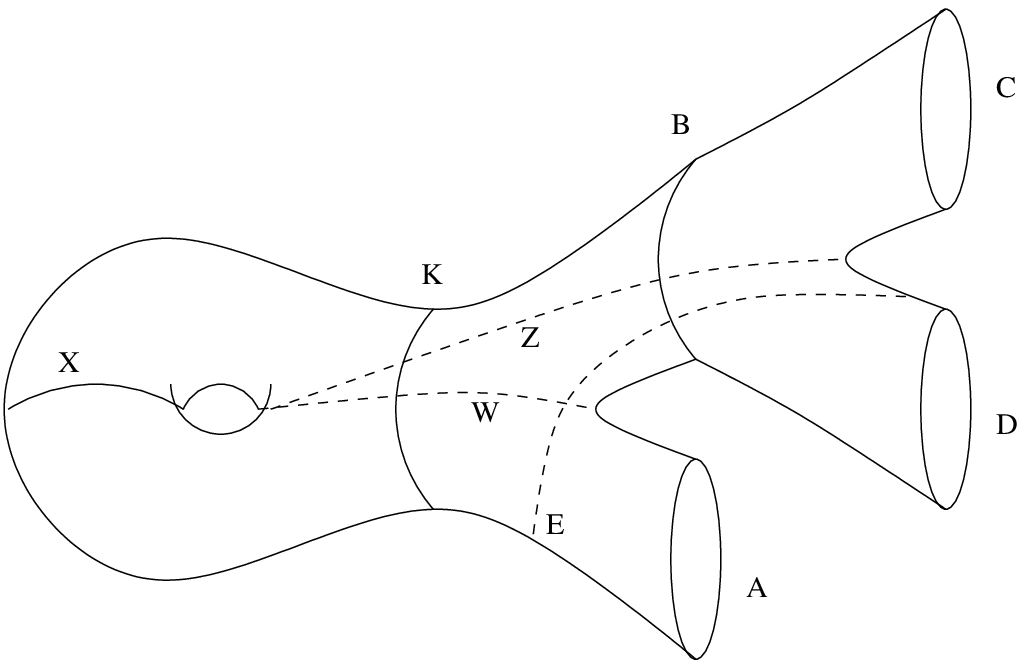}}
\centerline{F{\sc igure} 8: Three-holed torus.}
 
\

Suppose $B = \pm I.$
Then $C = \pm D^{-1}$ and
$A =\pm K.$
If $AD = E \neq \pm I$, then we may apply Proposition~\ref{prop:generic2hole}
to the two-holed torus
bounded by $E$ and $C$.  Hence,  if $E \neq \pm I,$ there exists an element in 
$\gamma \in \Gamma$ that fixes $a$ with $\gamma(k) \neq \pm a$
such that $\langle \gamma(X), \gamma(Y) \rangle$ generic.  
This implies $\gamma(B) \neq \pm I$. 
The same conclusion can be drawn if $AC \neq \pm I$.

Suppose that $AC= \pm I$ and  $AD = \pm I$.  Since  $A = \pm K,$
we have that 
$C,D = \pm K^{-1}.$
Hence, $B = CD= \pm K^{-2} = \pm I$.
This is only possible if $k=0$
or $K = \pm I.$ The latter case is ruled out by the generic
assumption on $\langle X, Y \rangle$.
Hence $k=0.$ We will show that in this special
case, we can obtain $\gamma(B) \neq \pm I$.

Suppose that $X$ commutes with $K$. If  $\tau_Y(X)= XY$ also commutes 
with $K$, then $X, Y, K$ all belong to the same one parameter subgroup
of $\SU$.  Hence $X$ and $Y$ commute. Since $\langle X, Y \rangle$ is generic,
we have that $XY$ does not
commute with $K.$
Suppose that $XKX^{-1}K^{-1} = -I$ ($X$ and $K$ anti-commute).  Then,
$X=K(-X)K^{-1},$ which implies that $x=0.$
Again, since $\langle X, Y \rangle$ is generic, we may assume that $x\neq 0.$
To summarize, since 
$\langle X, Y \rangle$ is generic,
one may arrange that $X$ neither commutes nor anti-commutes
with $K$ (i.e. $XKX^{-1}K^{-1}\neq \pm I$).

Now consider the curve $Z = XC.$
The Dehn twist in $Z$
preserves $E = \pm I$ and $c$.  Hence it fixes $k=0$
(consider the pants bounded by $C, E,$ and $K$).
On the four-holed sphere bounded by $C, D, X$ and $W$, the action
of $\tau_{Z}$ on the curve $B$ is 
$\tau_{Z}(B)= C(XC)D(XC)^{-1},$ so
$\tau_{Z}(B)= \pm K^{-1}XK^{-1}K^{-1}KX^{-1}= K^{-1}XKX^{-1} \neq \pm I.$
Hence $\gamma(B) \neq \pm I.$
Observe that by Remark \ref{rem:generic}
we may assume that $x$ is initially not in
$\{0, \pm 1, \pm \frac{1 \pm \sqrt{5}}{2}, \pm \sqrt 2\}$. 
Since $\tau_{Z}$ preserves $k$ and $x$,
$\langle \tau_{Z}(X), \tau_{Z}(Y) \rangle$
is also generic, note that
$\langle \tau_Z(X), \tau_Z(Y) \rangle$ is not $\Pin$ since $k=0$
and $x \not \in  \{\pm \sqrt 2, 0 \}.$
This proves the case $m=3$.
 
\

\centerline{\epsffile{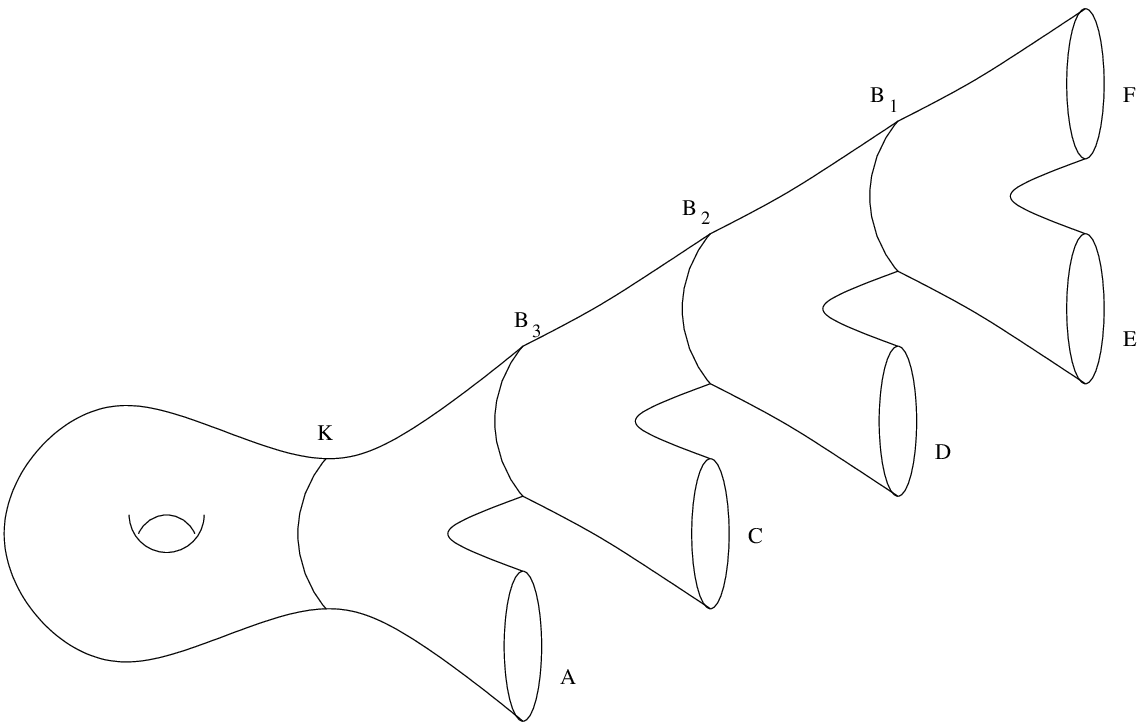}}
\centerline{F{\sc igure} 9: Getting rid of $\pm I$ on $\mathcal P$ for $m=5$.}

\

The above argument may be repeated iteratively, starting with loops in 
$\mathcal P$
that bound two boundary loops and working inward towards
$K.$  We demonstrate using the
case  $m=5$ with notation provided in Figure 9.

First consider the three-holed torus bounded by $E$, $F,$ and $ACD.$
Since $\langle X, Y \rangle$ is generic,
if $B_1= \pm I,$ then $ACD \neq \pm I$ (otherwise $K = \pm I$).
Use the previous argument for $m=3$ to arrange  
$B_1 \neq \pm I$.  
Next, the three-holed torus bounded by $B_1$, $D,$ and $AC$
is used to make $B_2 \neq \pm I$, etc.

\end{proof}

\begin{prop} \label{prop:generic3hole}
Suppose that $g=1,$ $m=3$ and $\rho$ is a generic representation.
Then the $\Gamma$-orbit $\Gamma([\rho])$ is dense in $\M_{\mathcal C}$.
\end{prop}

\begin{proof} Let $\epsilon' > 0$.
Suppose $[\rho], [\rho_0] \in \M_{\mathcal C}$ with
$[\rho]$ generic.  We adopt the usual practice of omitting
the symbol $\rho$ and using the subscript $0$ to denote
the values of $\rho_0$. We will also adopt the notation
of Figure 8. 
By Proposition  \ref{prop:dense}, without loss of generality
we may take $[\rho_0] \in \mathcal M_{\mathcal P}$ and show that
there exists $\gamma \in \Gamma$ such that
${\rm d} (\gamma([\rho]), [\rho_0])< \epsilon'.$

The map $f_{\mathcal P}$ is a submersion on $\M_{\mathcal P}$.
By the implicit function theorem, 
there exists  $\epsilon'>\epsilon>0,$ such that the
$\epsilon$-tubular neighborhood of $f_{\mathcal P}^{-1}(x_0,k_0,b_0) \subset 
\mathcal M_{\mathcal C}$ 
looks like
$[x_0-\epsilon, x_0+\epsilon] \times [k_0-\epsilon, k_0+\epsilon]$
$ \times [b_0-\epsilon, b_0+\epsilon]$
$ \times T^3.$

By Proposition \ref{prop:generic}, assume that $\langle X, Y \rangle$
is generic.
By
Proposition~\ref{prop:noI2}, we may assume $X,K,A,B,C,$ and $D$ 
are not
$\pm I.$  

%
%
%

Cutting along $X$, $K$ and $B$ yields a pants decomposition.  Hence,
by Proposition~\ref{prop:generic2hole} and 
Corollary~\ref{cor:densebase}, we only need to show that
there exists $\gamma \in \Gamma$ such that $\gamma([\rho])$
satisfies 
$|\gamma(b) - b_0|< \epsilon$,
with $\langle \gamma(X), \gamma(Y) \rangle$ generic.

The goal is to move $b$ near $b_0$.  The strategy is
to first move $x$ and $w$ near zero so that the moduli
space of the four-holed sphere 
bounded by $X,W,C,$ and $D$ contains $b$-coordinates
near $b_0$.  At the same time, one needs to ensure that
either $\tau_Z$ or $\tau_{Z''}$ has enough points on its orbit in
$E_{(c,d,x,w)}$, so that the $b$-coordinate
can actually be moved near $b_0$.  
This can be accomplished by
using the generic $\langle X, Y \rangle$ to
move the $y$-coordinate to a general position (i.e., outside 
of $Y_N$ for some large $N$).  This
will potentially produce large $\tau_Y$-orbits in $\langle Y, Z \rangle$ and 
$\langle Y, Z'' \rangle$.  Under such conditions, 
the $\tau_Y$-orbit in $\langle Y, Z \rangle$ and $\langle Y, Z'' \rangle$
will potentially contain points such that the $\tau_Z$- and 
$\tau_{Z''}$-orbits at these
points have a sufficiently large
number of points in $E_{(c,d,\gamma(x),\gamma(w))}$.  
For this strategy to work, however,
we must first deal with the special case where both
$\langle Y,Z \rangle$ and $\langle Y,Z''\rangle$ 
are $\Pin.$  These 
are the representations that are fixed by $\tau_Y$ (see 
Remark~\ref{rem:rotation1}).  In addition, one must also deal with
the analogous situation with the four-holed sphere (bounded by $X,W,C,$ and 
$D$) to ensure that either $\tau_Z$ or $\tau_{Z''}$ does not
fix the $b$-coordinate.  We begin with a detailed analysis of these 
exceptional cases where both of the following hold:
\begin{enumerate}
\item $\langle Y,Z \rangle$ is $\Pin$ or $\tau_Z$ fixes 
$(b, z, z'')$ 
\item $\langle Y,Z'' \rangle$ is $\Pin$ or $\tau_{Z''}$ fixes 
$(b, z, z'').$ 
\end{enumerate}

Case 1: Both $\langle Y,Z \rangle$ and $\langle Y,Z''\rangle$
are $\Pin$ and $y \neq 0$. This implies that
$z=z''=0$ and $\tr(YZ'')=\tr(YZ)=0.$
Note that $(b,z,z'') \in E_{(x,w,c,d)}.$ 
%
Note that the application of $\tau_B$ fixes the moduli space of
the two-holed torus bounded by $A$ and $B,$
i.e., $\tau_B$ fixes: $a$, $b$, $x$, $y,$ $k$, and $w.$
For the Dehn twist $\tau_B$ to fix $(b,z,z'')$ 
with $z=z''=0$ amounts to 
\begin {equation} \label{eq:c1}
wd+xc=0 \
\end {equation}
and
\begin {equation} \label{eq:c1.5}
xd+wc=0.
\end {equation}
In the special case where both $c=d=0,$ 
the defining equation for $E_{(x,w,c,d)}$
yields 
\begin {equation} \label{eq:c2}
b^2=xwb-x^2-w^2+4.
\end{equation}

Case 2: $\langle Y,Z'' \rangle$ is 
$\Pin,$ while $\langle Y,Z \rangle$ is not.
If $\tau_Z$ fixes $(b, z, z'') \in E_{(x,w,c,d)}$ then 
\begin{equation} \label{eq:c3}  
2b=xw+cd. 
\end{equation}

Case 3: $\langle Y,Z \rangle$ is 
$\Pin$ while $\langle Y,Z'' \rangle$ is not, with
$\tau_{Z''}$ fixing $(b, z, z''),$ then $2b=xw+cd$.
The argument for this case is symmetric to that of case 2.

Case 4: 
Neither $\langle Y,Z \rangle$ nor $\langle Y,Z'' \rangle$ is $\Pin$.
If both $\tau_Z$ and $\tau_{Z''}$
fix $(b,z,z'')$, then
$$
\left\{
\begin{array}{lll}
2b & = & xw+cd-zz''\\[2ex]
2z & = & xd+wc-bz'' \\[2ex]
2z'' & = & xc+wd-bz. 
\end{array}
\right.
$$
These lead to 
\begin{equation} \label{eq:c4}
2b-cd=wx - \frac {\bigl [ (dw+xc)-\frac b 2(xd+wc)\bigr ]
\bigl [(xd+wc)-\frac b 2 (dw+xc)\bigr ]
}{(2-\frac {b^2}2)^2}.
\end{equation}
A special case of the above equation is when $2b=cd.$ Then
equation (\ref{eq:c4}) becomes
\begin{equation} \label{eq:c5}
xw= \frac {\bigl [ (dw+xc)-\frac b 2(xd+wc)\bigr ]
\bigl [(xd+wc)-\frac b 2 (dw+xc)\bigr ]
}{(2-\frac {b^2}2)^2}.
\end{equation}
Note this quadratic equation is degenerate in $w$ 
only if $c=\pm 2,$ $d=\pm 2,$ $c=b=0,$ or $d=b=0.$ 
The case  $c=b=0$
leads to the equations $2z=xd,$ $2z''=wd,$ and
$zz''=xw.$
Therefore,  $4zz''=xwd^2,$  i.e., $d^2=4.$
Similarly, $d=b=0$ leads to $c^2=4.$
Note that the non-degenerate solutions of equation (\ref{eq:c5})
together with $2b=cd$
are: 
\begin{equation} \label{eq:c6}
w= cdx/4 
\end{equation}
and
\begin{equation} \label{eq:c6.5}
x= cdw/4.
\end{equation}

Now we carefully choose a generic representation 
of the two-holed torus bounded by $A$ and $B$
that allows us to avoid the special cases presented above.
\begin{lemma} \label{lem:laundry2}
For any integer $J>0$, 
there is
$\gamma \in \Gamma$ such that $\gamma (\rho)$ has
$\tau_Y$-orbit 
containing at least $J$ points satisfying
all of the following:
\begin{enumerate}
\item Each of the $J$ points have distinct $x$ and
$w$-coordinates inside $(-\delta  ,\delta)$, where
$\delta$ is given as in Lemma~\ref{lem:inband} applied to the
four-holed sphere that bounds $X$, $W$, $C$, and $D$.
\item $|x|, |w|$ are small so that 
$b_0 \pm \frac \epsilon 2 \in I_{c,d}\cap I_{x,w}.$
\item If $c \neq 0$ or $d \neq 0$, then these $J$ points
do not belong to the subvariety defined by either equation~(\ref{eq:c1})
or (\ref{eq:c1.5}).
\item If $c=d=0$, then these $J$ points do not belong to the 
subvariety defined by equation~(\ref{eq:c2})
\item The $J$ points do not belong to the subvariety
defined by equation~(\ref{eq:c3}).
\item The $J$ points do not belong to the subvariety
defined by equation~(\ref{eq:c4}) if $2b \neq cd$ and
do not belong to the subvarieties defined by equations~(\ref{eq:c6})
and (\ref{eq:c6.5})
if $2b = cd$.
\end{enumerate}
Note that the $J$ points provided in Lemma~\ref{lem:laundry2} are obtained by 
using the mapping class action on the two-holed torus bounded by $A$ and $B$.  
Hence the value $b \neq \pm 2$ remains fixed.
\end{lemma}
\begin{proof}
We begin by noting that the two-holed torus
bounded by $A$ and $B$ is generic.

Conditions 1 and 2 follows 
from Proposition~\ref{prop:generic2hole}
and Remark ~\ref{rem:generic}.
Note that both $|x|$ and $|w|$ can be made simultaneously near zero,
since $x$ can always be made arbitrarily close
to zero for any value of $k$ (see Remark ~\ref{rem:generic})
and since the $w$ values of
$E_{\kappa}$ for $\kappa=(a,b,x,x)$ 
take on
values inside $I_{x,a}\cap I_{x,b},$ which, for $x$ near zero, 
contains the value $w=0.$

By Remark~\ref{rem:generic2}, for $y$ having an $\eta$-dense orbit for $\eta$
sufficiently small,
we may find a point with a $\tau_Y$-orbit with at least $J$ points
not belonging to any of the proper subvarieties described in 
conditions 3-6. Moreover, after a possible application of $\tau_B,$ we
also have that one of $\langle Y, Z \rangle$ or 
$\langle Y, Z''\rangle$ is not $\Pin$.
\end{proof}

%
%

Conditions 3-6 of Lemma~\ref{lem:laundry2} ensure that 
for every point on the $\tau_Y$-orbit of $\gamma(\rho)$
either
$\langle Y, Z \rangle$ is not $\Pin$ and $\tau_Z$ does not fix $(b, z, z'')$ 
or 
$\langle Y,Z'' \rangle$ is not $\Pin$ and $\tau_{Z''}$ does not fix
 $(b, z, z'').$ 
%
%
%
We assume without loss of generality assume the former.

Finally, we start with a point $\gamma(\rho)$ with
$J$ sufficiently large so that one of the
points on the $\tau_Y$-orbit of $[\gamma(\rho)]$ has a $\tau_Z$-action
containing  at least 
$$N(\frac {\delta}{N( \frac \epsilon 2)+3})+3$$ points.

The argument now follows similarly to that presented in 
Proposition~\ref{prop:generic2hole} to produce $\gamma \in \Gamma$
by obtaining $|\gamma(b)-b_0| < \epsilon.$ 
This completes the proof.

To summarize, one first makes $b \neq \pm 2$.  Then uses the
fact that the two-holed torus bounded by $A,B$ is generic to move
the $x,w$-coordinates near zero.  This allows the $b$-coordinate
to be moved near $b_0$.  Finally the generiticity of the two-holed
torus allows one to move the $x,w$-coordinates near $x_0,w_0$.
\end{proof}

\begin{prop} \label{prop:genericnhole}
Suppose $g=1$ and $M$ has $m \ge 1$ boundary components
with $\rho$ a generic representation.
Then the $\Gamma$-orbit $\Gamma([\rho])$ is dense in $\M_{\mathcal C}$.
\end{prop}
\begin{proof}
The cases $m=1, 2,$ and $3$ have previously been established.
For $m \ge 4,$ the argument
follows by an induction process similar to that
in the proof of Proposition~\ref{prop:noI2}.

Let $\epsilon '> 0$.
Let  $\rho$ be a generic representation. By
Proposition \ref{prop:generic}, $M$ has a 
generic handle $\langle X, Y \rangle$

We demonstrate how to proceed in the case $m=5$
(see the Figure 10 below).
By Proposition  \ref{prop:dense}, without loss of generality
we may take $[\rho_0] \in \mathcal M_{\mathcal P}$ and show that
there exists $\gamma \in \Gamma$ such that
${\rm d} (\gamma([\rho]), [\rho_0])< \epsilon'.$

The map $f_{\mathcal P}$ is a submersion on $\M_{\mathcal P}$.
By the implicit function theorem, 
there exists $\epsilon'>\epsilon>0,$ such that the
$\epsilon$-tubular neighborhood of $$f_{\mathcal P}^{-1}\bigl ((x_0,k_0,(b_1)_0,(b_2)_0,(b_3)_0)\bigr ) \subset 
\mathcal M_{\mathcal C}$$ 
looks like
$B_{\epsilon}(\{(x_0,k_0,(b_1)_0,(b_2)_0,(b_3)_0\})) \times T^5,$
where $B_{\epsilon}(S)$ denotes the  $2\epsilon$-open box 
neighborhood of the set $S$.

\

\centerline{\epsffile{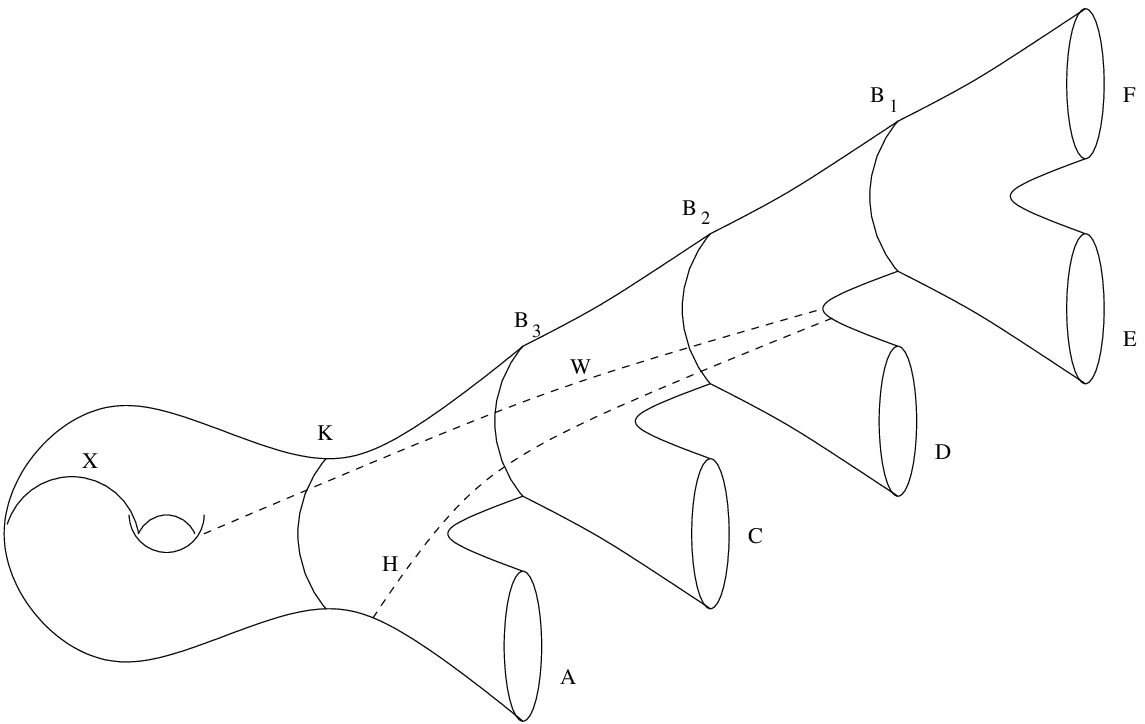}}
\centerline{F{\sc igure} 10: The induction on the 5-holed torus.}

\


By Proposition~\ref{prop:noI2}, we have that 
$B_1 \neq \pm I.$  By induction, Proposition~\ref{prop:genericnhole}
is true for the four-holed
torus bounded by $B_1, A, C$ and $D$.
We use this
four-holed torus
to arrange for $x$ and $w$ (see Figure 10) 
to have near zero traces. This ensures that 
the $b_1$-coordinate of $[\rho_0]$ is accessible. 
We then arrange $H \neq \pm I.$
We now cut $M$ at $H$ and use Proposition~\ref{prop:generic3hole}
to get the $b_1$-coordinate of $\gamma([\rho])$ within
$\epsilon$ of $(b_1)_0.$
Next, we get the $b_2$-coordinate  of $\gamma([\rho])$ within
$\epsilon$ of $(b_2)_0,$ by using the
four-holed torus obtained by cutting at $B_1$.
Finally, Proposition ~\ref{prop:generic3hole} and
Corollary~\ref{cor:densebase} yields the result.
The situation is identical for any $m>3.$
\end{proof}

\section{Genus $g$ with $m$ boundary components}\label{sec:induction}

Suppose $M$ is a surface with $g>1$ and $m \ge 0.$  Again,
we assume that $C_i \neq \pm I$ for all boundary curves $C_i \in \partial M,$
unless $m=1,$ $g>1,$ and $C_1 = \pm I.$

\begin{prop} \label{prop:noI}
Let  $M$ have
generic  $\langle X, Y \rangle$ that is part of a 
pants decomposition $\mathcal P.$ Then
there is $\gamma \in \Gamma$
 so that
$\gamma (\rho(P)) \neq \pm I,$ for all $P \in \mathcal P$
with $\langle \gamma(X), \gamma(Y) \rangle$ generic.
If $m=1$ and $C_1=\pm I,$ then 
$\gamma (\rho(P)) \neq \pm I$ for all $P \in \mathcal P$
except $C_1.$ Furthermore, 
$\langle \gamma (P),\gamma(Q)\rangle$ is generic, 
for all $P \in \mathcal P$ 
and $PQP^{-1}Q^{-1}=R \in \mathcal P$ that bound a
three-holed sphere that forms a one-holed torus in $M.$ 
\end{prop}

\

\centerline{\epsffile{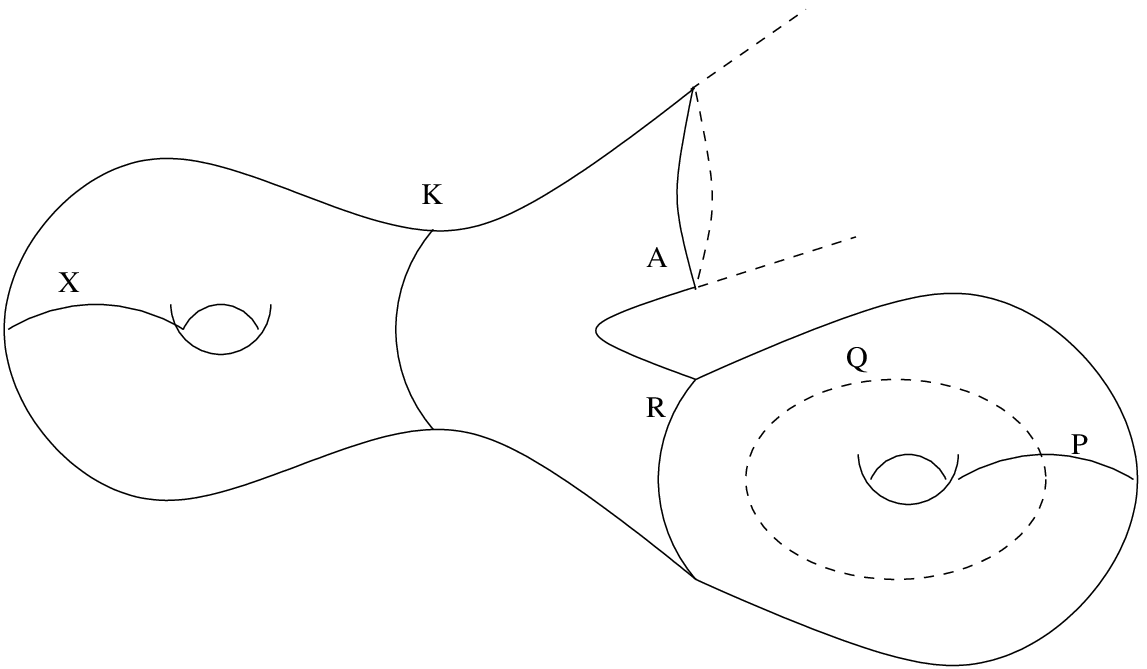}}
\centerline{F{\sc igure} 11: A punctured genus two inside $M$.}

\

\begin{proof}
We first show that there is $\gamma \in \Gamma$ so that
$\langle \gamma (P),\gamma(Q)\rangle$ is generic, 
for all $P \in \mathcal P$ 
and $PQP^{-1}Q^{-1}=R \in \mathcal P$ that bound a
three-holed sphere that forms a one-holed torus in $M.$ 
Consider $P,$ $Q,$ and $R$ as in Figure 11, 
where $P, R \in \mathcal P$ bound a
pants that forms a one-holed torus in $M.$

Case 1: Suppose that $R \neq \pm I.$ 
In this case, at least one of $p$, $q$ 
or $\tr(PQ)$
is not in $\{0, \pm 2\}.$ 
Without loss of generality, assume that 
 $p \notin \{0, \pm 2\} $ (otherwise, Dehn twist
in $Q$). Thus, we may cut $M$ at $P$ and $A$
and apply Proposition~\ref{prop:generic3hole} to the resulting
three-holed torus to obtain $r \not\in \mathcal S$ with
$r \neq p^2+2.$ Note that in the special case $A  = \pm I,$ we  apply
Proposition~\ref{prop:generic2hole} along with Remark \ref{rem:reduction}.
Thus the handle $\langle \gamma (P), \gamma(Q) \rangle$ is generic. 

Case 2: Suppose that $R = - I.$ Then both $P$ and $Q$ are not $\pm I.$
Moreover, $A \neq \pm I$ since  $K \neq \pm I.$ Thus, we may cut $M$ at $P$
and $A,$ and then apply Proposition~\ref{prop:generic3hole}
to reduce the situation to Case 1.

Case 3:  Suppose that $R = I.$ Then, necessarily, $A = \pm K \neq \pm I$ since
$\langle X, Y \rangle$ is generic. 
If either $P\neq \pm I$ or $Q \neq \pm I,$ we may cut at $P$
(or $Q$) and apply Proposition~\ref{prop:generic3hole}
as in case 2.
Since  $PA= \pm K \neq \pm I$,
If both  $P = \pm I$ and $Q = \pm I,$ then
$\tau_{PA}(Q) = Q(PA) = \pm  K \neq \pm I.$
Note that the action of $\tau_{PA}$ does not affect the
generic handle $\langle X, Y \rangle.$
We now cut at
$Q$ and apply Proposition~\ref{prop:generic3hole}
as in case 2.
This argument can be applied independently to each of the 
$g-1$ handles of $M$.

Having obtained that all handles are generic, we now apply
Proposition~\ref{prop:noI2} to the remaining curves
in $\mathcal P$ that are interior to the $(m+g-1)$-holed torus
obtained by cutting off each of the $g-1$ handles.
In the special case $m=1$ and $C_1=\pm I,$  
 Proposition~\ref{prop:noI2} applies to those curves of $\mathcal P$
that are separated from $C_1$ by the curve $KC_1.$ Moreover, 
$KC_1 \neq \pm I$ since $\langle X, Y \rangle$ is assumed generic.
\end{proof}


\

\centerline{\epsffile{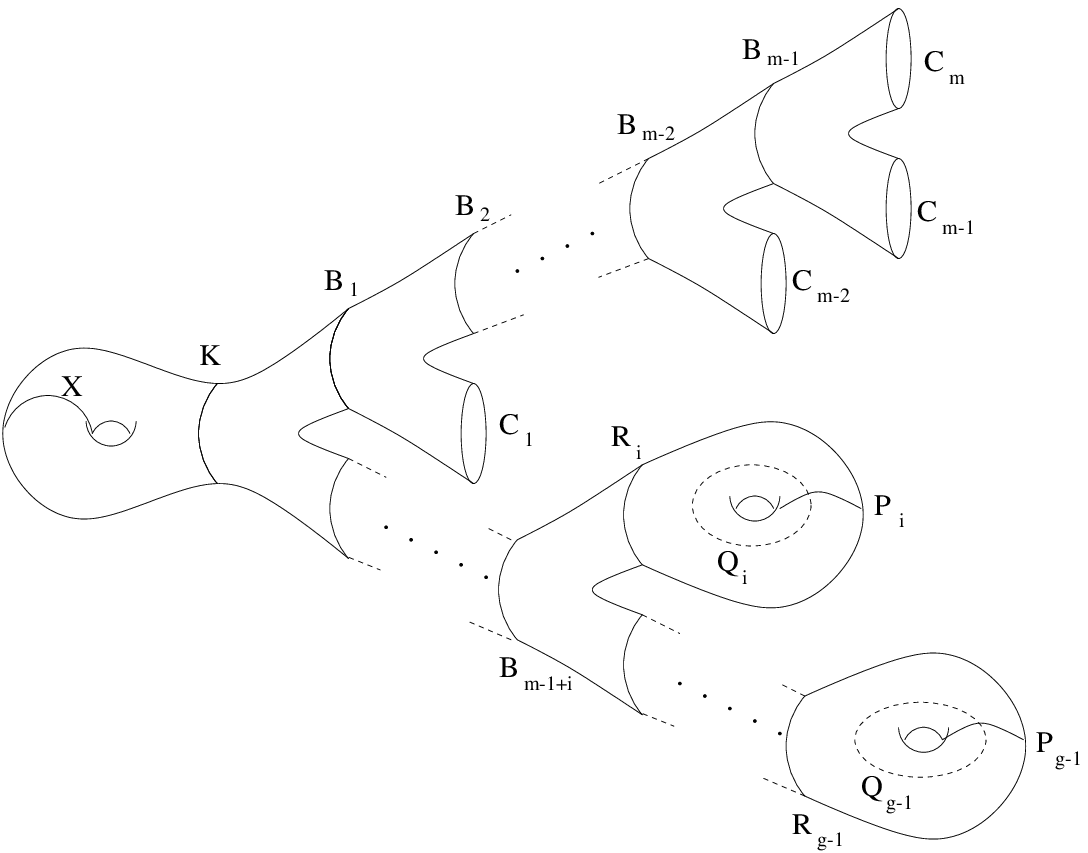}}
\centerline{F{\sc igure} 12:  Decomposition of $M$.}

\

Now we prove Theorem~\ref{thm:3}.
We adopt the notation of Figure 12 for the decomposition of $M$.
Let $\epsilon' > 0$.
Let  $\rho$ be a generic representation. By
Proposition \ref{prop:generic}, $M$ has a 
generic handle $\langle X, Y \rangle$
(we adopt the notation presented in Figure 12).
By Proposition  \ref{prop:dense}, without loss of generality
we may take $[\rho_0] \in \mathcal M_{\mathcal P}$ and show that
there exists $\gamma \in \Gamma$ such that
${\rm d} (\gamma([\rho]), [\rho_0])< \epsilon'.$

The map $f_{\mathcal P}$ is a submersion on $\M_{\mathcal P}$.
By the implicit function theorem, 
there exists  $\epsilon'>\epsilon>0,$ such that the
$\epsilon$-tubular neighborhood of $f_{\mathcal P}^{-1}(\beta) \subset 
\mathcal M_{\mathcal C},$ where
$$\beta= (x_0,k_0,(r_1)_0,...,(r_{g-1})_0,(p_1)_0,...,(p_{g-1})_0,(b_1)_0, ..., (b_{(m-1)+(g-2)})_0 ),$$
looks like
$$B_{\epsilon}(\{\beta\})    \times T^{3g+m-3},$$
where $B_{\epsilon}(S)$ denotes the  $\epsilon$-open box 
neighborhood of the set $S$.

%

By Proposition~\ref{prop:noI}, we may assume that 
$p_i \not \in \{0, \pm 2 \}$ for all $1\le i \le g-1$. 
Cutting $M$ along all of the $P_i$'s results in a
$(m+2g-2)$-holed torus
(In the case of
$m=1,$  $C_i= \pm I,$ and $g>1$, 
Remark~\ref{rem:reduction} and Proposition~\ref{prop:genericnhole}
give us a generic $(2g-2)$-holed torus).

Dehn twist each generic handle $\langle P_i, Q_i \rangle$ 
to obtain $p_i$
arbitrarily close to zero. Since $x$ and
$z_i=\tr(R_iX)$ can be made arbitrarily close to zero by
Proposition~\ref{prop:genericnhole}, we
are ensured that the target value $r_i$-coordinate of $\rho_0$ 
will be  inside
the $\epsilon-$neighborhood of
$I_{p_i,p_i} \cap I_{x,z_i}.$ Thus, there are points on the four-holed sphere
$E_{(p_i,p_i,x,z_i)}$ with $r_i$-coordinates within
$\epsilon$ of $(r_i)_0$, the $r_i$-coordinate of $\rho_0.$ 

We next re-apply Proposition~\ref{prop:genericnhole} to the 
$(m+2g-2)$-holed torus to move each $r_i$ 
within $\epsilon$ of $(r_i)_0,$ keeping $\langle \gamma(P_i), \gamma(Q_i) 
\rangle$ generic
by ensuring that $\gamma(p_i) \neq 0,$ $\gamma(r_i) \not\in {\mathcal S}$ and 
$\gamma(r_i) \neq (\gamma(p_i))^2-2$.
We start working outward in (as suggested in Figure 12).
Next, we  twist on each of the $g-1$ generic handles to get $p_i$ within
$\epsilon$ of  $(p_i)_0.$
Finally, we cut off each handle at $R_i$
and apply Proposition~\ref{prop:genericnhole} to the resulting 
$(m+g-1)$-holed torus with generic 
$\langle X, Y \rangle$. 

To summarize, we first moves the $p_i,r_i$-coordinates 
away from $\pm 2$.
Then, we cut along the 
$P_i$'s and use the resulting generic $(m+2g-2)$-holed torus
to move each $r_i$-coordinate to make each $\langle P_i, Q_i \rangle$ generic,
while preserving the generiticity of $\langle X, Y \rangle$.  Next,
we use the generic torus $\langle P_i, Q_i \rangle$ to
move the $p_i$-coordinates near zero.  As $x$ and $z_i=\tr(R_iX)$
can be taken
arbitrarily close to zero as well, we see that
there are points on the four-holed sphere
$E_{(p_i,p_i,z_i,x)}$ with $r_i$-coordinate within
$\epsilon$ of $(r_i)_0.$
Again by using the 
generic $(m+2g-2)$-holed torus, we may move the $r_i$-coordinate
near $(r_i)_0$.  After that,
we move the $p_i$-coordinate near $(p_i)_0$.
Finally cut at the $R_i$'s and use the generic $(m+g-1)$-holed torus
to move the $b_i$-coordinates near $(b_i)_0$
and lastly, the $x,k$-coordinates near $x_0,k_0$.
Theorem~\ref{thm:3} now follows from Corollary~\ref{cor:densebase}.


\appendix
\section{Proof of Proposition~\ref{prop:generic}}
Much of the proof is computational,  repetitive, and performed 
with the aid of MAPLE. However, each sub-case can be worked out
by hand. 
In what follows, we provide several typical examples for the various cases;  
however the MAPLE files for all cases are available online at
 \underline{http://vortex.bd.psu.edu/$\sim$jpp/td2/}.

We begin by defining two constants:
$$
\left\{
\begin{array}{lll}
r & = & \frac{\sqrt{5} + 1}{4} \\[2ex]
s & = & \frac{\sqrt{5} - 1}{4}.
\end{array}
\right.
$$
Suppose that $\rho \subset \Hom(\pi_1(M,O), \SU)$ is generic and
$\langle A_i, A_{g+i} \rangle$ is $G$ for some
proper subgroup $G \subset \SU$.  Suppose that there exists an element
$A_j \in \SU$ such that $A_j \not\in G$.  If
$O$ is chosen to be the intersection point of $A_i$ and $A_{g+i}$,
then for any $j$, one may construct a loop corresponding to 
the fundamental group element
$A_{g+i} A_j$ such that the loops $(A_i, A_{g+i} A_j)$ 
intersect only at $O$.  There is an analogous construction
for $(A_i, A_jA_{g+i}),$  as well as for
$(A_i A_j, A_{g+i})$ and
$(A_j A_i, A_{g+i})$.

\begin{lemma}
There exists a handle $(A_i,A_{g+i})$ such that 
either $A_i \neq \pm I$ or $A_{g+i} \neq \pm I.$
\end{lemma}
\begin{proof}
Suppose that $A_i = \pm I$ for all $1 \le i \le 2g.$ 
Since $\rho$ is generic, there exists $A_{j}\neq \pm I$ for
$j>2g.$  
The handle  $(A_i A_j, A_{g+i})$ has the desired property.
\end{proof}

\begin{lemma}
There exists a handle $(A_i,A_{g+i})$ such that 
$\rho\vert_{\langle A_i, A_{g+i} \rangle}$ is not $\Spin$.
\end{lemma}
\begin{proof}
Since $\rho$ is generic, by the previous Lemma, there
exists $A_i \neq \pm I$ for $1 \le i \le 2g.$ 
Then $A_i$ defines a one-parameter subgroup $P$ of $\SU$.  
Move $O$ to the intersection point of $A_i$ and $A_{g+i}$.  Since $\rho$
is generic, there exists $A_j \not\in P$.  This implies that
$\langle A_i, A_{g+i} A_j \rangle$ is non-Abelian, hence, not $\Spin$.
\end{proof}

\begin{lemma} \label{lem:pin}
There exists a handle $(A_i,A_{g+i})$ such that 
$\rho\vert_{\langle A_i, A_{g+i} \rangle}$ is not $\Pin$.
\end{lemma}
\begin{proof}
Suppose $\langle A_i, A_{g+i} \rangle$ is $\Pin$, but not $\Spin$. 
The group $\Pin$ acts on $S^2$ 
via its quotient ${\rm SO}(3)$ and 
preserves a circle $S^1 \subset S^2$.  There are two cases.

\noindent{\sl Case 1:} $\langle A_i, A_{g+i} \rangle$ preserves
a unique circle $S^1 \subset S^2$. 
Since $\rho$ is generic, there exists
$A_j$ not preserving $S^1$.  This implies that 
$\langle A_i, A_{g+i} A_j \rangle$ does not preserve
any  circle in $S^2$, hence, not $\Pin$.

\noindent{\sl Case 2:} $\langle A_i, A_{g+i} \rangle$ preserves
three circles: $S_1, S_2, S_3$.  Then $\tr(A_i) = 
\tr(A_{g+i})=0$ and  $\tr(A_iA_{g+i}) =  0.$ 
This implies that $A_i^2=-I$, $A_{g+i}^2=-I,$ and
$A_iA_{g+i} = -A_{g+i}A_i.$

Suppose further that $\langle A_i, A_{g+i} A_j \rangle$ and
$\langle A_iA_j, A_{g+i}\rangle$ are
$\Spin$ or $\Pin$ preserving
three circles.  Then
up to conjugacy,
$$ A_i = \left[
\begin{array}{cc}
i  & 0\\
0& -i\\
\end{array}
\right],\ \ 
 A_{g+i}= \left[
\begin{array}{cc}
0  & 1 \\
-1 & 0\\
\end{array}
\right],\ \ \text{\rm  and }  
A_iA_{g+i} = \left[
\begin{array}{cc}
0 & i\\
i & 0\\
\end{array}
\right].\ \ 
$$ 
Let $$A_j =  \left[
\begin{array}{cc}
w + ix  & y + iz\\
-y+ iz &  w - ix\\
\end{array}
\right],\ \ 
$$
with $w^2+ x^2 + y^2 + z^2 = 1$ be such that 
$A_j$ is not contained in the isomorphic copy of 
$\langle A_i, A_{g+i} \rangle$. 

Since  $\langle A_i, A_{g+i} A_j \rangle$ is either 
$\Spin$ or $\Pin$ preserving
three circles, we have that $A_iA_{g+i} A_j = \pm A_{g+i} A_j A_i.$
However, this implies that  $-A_{g+i}A_i A_j = \pm A_{g+i} A_j A_i,$ so
$A_i A_j = \pm A_j A_i.$ For 
$A_i A_j = A_j A_i,$ we see that  $y=z=0.$
For 
$A_i A_j = -A_j A_i,$ we see that  $x=w=0.$

Likewise, since $\langle A_iA_j, A_{g+i}\rangle$ is either $\Spin$ 
or $\Pin$ preserving
three circles, we have that $A_{g+i} A_j = \pm A_j A_{g+i}.$
The case
 $A_{g+i} A_j = A_j A_{g+i}$ leads to $x=z=0,$ while
the case $A_{g+i} A_j = -A_j A_{g+i}$ leads to $w=y=0.$

Taken together, we see that in all cases three of 
$w, x, y, z$ are zero. This implies that 
$A_j$ is equal to one of: $\pm A_i, \pm A_{g+i},$  
$\pm A_iA_{g+i},$ or $\pm I,$ which is a 
contradiction. 
\end{proof}

\begin{lemma} \label{lem:tetra}
There exists a handle $(A_i,A_{g+i})$ such that $\rho\vert_{\langle A_i, A_{g+i} \rangle}$
is neither $T$ nor $\Pin$.
\end{lemma}
\begin{proof}
By Lemma~\ref{lem:pin}, there exists $i$ such that
$\langle A_i, A_{g+i} \rangle$ is not $\Pin$.
Suppose $\langle A_i, A_{g+i} \rangle$ is $T$.
By applying Dehn twists, and in light of the analysis carried
out in \cite{Pr1},
one may assume that  $\rho$ has character
$(\tr(A_i),\tr(A_{g+i}),\tr(A_iA_{g+i})) =(1, \pm 1, \pm 1)$
on $(A_i,A_{g+i}),$ with signs taken together.
We will only address representations with character
$(1,1,1),$ as the arguments for $(1,-1,-1)$ are practically identical.

Up to conjugation, we may choose:
$$A_i = \frac 1 2\left[
\begin{array}{cc}
1 - i  & 1 - i\\
-1 - i & 1 + i\\
\end{array}
\right], \ \ 
A_{g+i} = \frac 1 2\left[
\begin{array}{cc}
1 - i  & -1 + i\\
1 + i & 1 + i\\
\end{array}
\right].\ \ 
$$

The group $T=\langle A_i, A_{g+i} \rangle$
consists of: 
$\pm I,$ the 16 matrices of the form
$$ \pm \frac 1 2\left[
\begin{array}{cc}
1 \pm i  & \pm 1 \pm i\\
\mp 1 \pm i & 1 \mp i\\
\end{array}
\right],\ \ 
$$ 
(signs in the second row determined by those in the first),
and the 6 matrices of the form
$$ \pm \left[
\begin{array}{cc}
i  & 0\\
0& -i\\
\end{array}
\right],\ \ 
\pm \left[
\begin{array}{cc}
0  & 1 \\
-1 & 0\\
\end{array}
\right],\ \  
\pm  \left[
\begin{array}{cc}
0 & i\\
i & 0\\
\end{array}
\right].\ \ 
$$ 
Let $$A_j =  \left[
\begin{array}{cc}
w + ix  & y + iz\\
-y+ iz &  w - ix\\
\end{array}
\right],\ \ 
$$
with $w^2+ x^2 + y^2 + z^2 = 1$ be such that 
$A_j$ is not contained in the isomorphic copy of 
$T = \langle A_i, A_{g+i} \rangle$. 

First, consider the case where each of
$\langle A_i, A_{g+i} A_j \rangle,$ 
$\langle A_i, A_jA_{g+i} \rangle,$  
$\langle A_iA_j, A_{g+i}\rangle,$ and
$\langle A_jA_i, A_{g+i} \rangle$  
are $\Spin$.
This means that $A_i(A_{g+i} A_j)=(A_{g+i} A_j)A_i
=(A_iA_j) A_{g+i},$ which implies that
$A_jA_{g+i}= A_{g+i} A_j.$
Similarly,
$A_jA_{i}= A_{i} A_j.$
This means that $A_j \in \langle A_i, A_{g+i}\rangle,$
which is a contradiction. 


Next, suppose that $\langle A_i, A_{g+i} A_j \rangle$ is  $\Pin$
and that $\langle A_i A_j, A_{g+i} \rangle$ is $\Spin$.
Then $\tr(A_{g+i}A_j) = w + x + y - z = 0$ and
$\tr( A_iA_{g+i} A_j)= w + x -y - z = 0$
implies that $y = 0$ and $z = w + x.$
The equation
$$\tr(A_iA_j)^2 + \tr(A_{g+i})^2 +\tr( A_i A_jA_{g+i})^2 - 
\tr(A_iA_j)\tr(A_{g+i})\tr( A_i A_jA_{g+i})-2 = 2,$$
implies that $z = \pm \frac{\sqrt 3} 2.$ 
Finally, $w^2 + x^2 + z^2 = 1 $ implies that 
either $x$ or $w$ is not real, a contradiction.
A similar argument holds in the event that
$\langle A_i, A_{g+i} A_j \rangle$ is $\Spin$
and $\langle A_i A_j, A_{g+i} \rangle$ is   
$\Pin.$

We indicate how to proceed in the case where 
$\langle A_i, A_{g+i} A_j \rangle$ is $T$ 
and $\langle A_i A_j, A_{g+i} \rangle$ is $\Spin$.  
Since $\langle A_i, A_{g+i} A_j \rangle$ is  $T,$ both
$\tr(A_{g+i} A_j)$ and
$\tr(A_iA_{g+i}A_j)$ must take on the values $\pm 1 $ or $0.$
The case where both are zero corresponds to a $\Pin$ representation
which was already handled.

Suppose that
$\tr(A_{g+i} A_j)=0$ and $\tr(A_iA_{g+i}A_j)=1.$
Then
$w+x+y-z = 0$ and $w+x-y-z = 1,$
so $y = -\frac{1}{2}$
and $z = w + x - \frac{1}{2}.$
The equation 
$$\tr(A_i A_j)^2 + \tr(A_i A_j A_{g+i})^2 -\tr(A_i A_j)\tr(A_i A_j A_{g+i})-3 = 0$$
yields $z=-1$ or $\frac{1}{2}$.
If $z=-1,$ then $w$ is not real.  This implies that
$w=\frac 1 2,$ $x = \frac 1 2,$ $y = -\frac 1 2,$ and
$z = \frac 1 2.$ However, this implies that 
$A_j \in \langle A_i, A_{g+i} \rangle$ which is a contradiction.
Similar arguments hold in the other 7 cases, as well as the 8 cases where
$\langle A_i, A_{g+i} A_j \rangle$ is $\Spin$
and $\langle A_i A_j, A_{g+i} \rangle$ is $T$.

We now address the remaining cases where
$\langle A_i, A_{g+i} A_j \rangle$ and $\langle A_i A_j, A_{g+i} \rangle$ 
are either $T$ or $\Pin$.
In these cases, 
$\tr(A_{g+i} A_j)$, $\tr(A_{i} A_j)$, and
$\tr(A_iA_{g+i} A_j)$ must take on the values $\pm 1 $ or $0.$

We explicitly work
out a $T$-$\Pin$ case in which $\tr(A_{i} A_j),$ $\tr(A_iA_{g+i}A_j),$ and
$\tr(A_{g+i} A_j)$ are all zero.
These lead to the equations
$$w+x-y-z=0$$
$$w+x-y+z=0$$
$$w+x+y-z=0$$
$$w^2+x^2+y^2+z^2=1.$$

Therefore, $y=0,$ $z=0$ and  $w=-x,$ which implies that $x = \pm \frac {\sqrt 2}2,$
and $w = \mp  \frac {\sqrt 2}2.$
However, the trace of the matrix $A_iA_{g+i}A_iA_j$ is $\pm \sqrt 2.$
Thus the handle
$\langle A_i, A_{g+i}A_iA_j \rangle$ is neither $T$ nor $\Pin$.
(This handle is obtained by first twisting in $A_i$ on 
$\langle A_i,A_{g+i}\rangle$, then 
using the loop $A_j.$)




The remaining 26 cases are obtained by varying the
values of  $\tr(A_{i} A_j),$ $\tr(A_iA_{g+i}A_j),$ and
$\tr(A_{g+i} A_j)$ in $\{0, \pm 1\}.$ The arguments proceed similarly
to those exemplified above and
lead to: 
\begin{enumerate}
\item matrices $A_j$ inside
$\langle A_i, A_{g+i} \rangle$ (a contradiction),
\item a handle $\langle A_i, A_{g+i}A_iA_j \rangle$
that is neither $T$ nor $\Pin,$
\item one or more of $w, x, y, z$ is not real (a contradiction).
\end{enumerate}
\end{proof}

\begin{lemma}
There exists a handle $(A_i,A_{g+i})$ such that $\rho\vert_{\langle A_i, A_{g+i} \rangle}$
is neither $C$ nor $\Pin$.
\end{lemma}
\begin{proof}
Assume that $\langle A_i, A_{g+i} \rangle$ is $C$ but neither 
$T$ nor $\Pin.$
By applying Dehn twists and in light of the analysis carried
out in \cite{Pr1},
one may assume that  $\rho$ has character
$(\tr(A_i),\tr(A_{g+i}),\tr(A_iA_{g+i})) =(1, \pm 1, \pm 1)$
on $(A_i,A_{g+i}),$ with signs taken together.
We will only address representations with character
$(\sqrt 2,1,\sqrt 2),$ as the arguments for $(\sqrt 2,-1,-\sqrt 2)$ are practically identical.

Up to conjugation, we may choose
$$A_i = \frac 1 2\left[
\begin{array}{cc}
\sqrt 2+ \sqrt 2 i  & 0\\
0 & \sqrt 2- \sqrt 2 i \\
\end{array}
\right],\ \ 
A_{g+i} = \frac 1 2\left[
\begin{array}{cc}
1 - i  & -1 - i\\
1 - i  & 1 + i\\
\end{array}
\right]. \ \ 
$$

The group $C = \langle A_i, A_{g+i} \rangle$ 
consists of: the 24 tetrahedral matrices in the previous lemma,
the four matrices 
$$ \pm \frac 1 2\left[
\begin{array}{cc}
\sqrt 2 \pm \sqrt 2i  & 0\\
0 & \sqrt 2 \mp \sqrt 2i\\
\end{array}
\right]\ \ 
$$ 
(signs in the second row determined by those in the first),
and the 20 matrices obtained by permuting
two non-zero terms ($\sqrt 2$) and
two zero terms ($0+0i$) in the first row of the matrix.
One such example is:
$$ \frac 1 2\left[
\begin{array}{cc}
\sqrt 2i & \sqrt 2 \\
-\sqrt 2 & -\sqrt 2i \\
\end{array}
\right].\ \ 
$$ 

Let $$A_j =  \left[
\begin{array}{cc}
w + ix  & y + iz\\
-y+ iz &  w - ix\\
\end{array}
\right],\ \ 
$$
with $w^2+ x^2 + y^2 + z^2 = 1$ be such that $A_j$ is
not contained in the isomorphic copy of $C$
generated by $\langle A_i, A_{g+i} \rangle.$
Note that
$\tr(A_iA_j) = \sqrt 2 w - \sqrt 2 x,$
$\tr(A_{g+i}A_j) = w +x +y + z,$
$\tr(A_iA_jA_{g+i}) = \sqrt 2 w + \sqrt 2 y,$
 and
$\tr(A_iA_{g+i}A_j) = \sqrt 2 w + \sqrt 2 z.$

If each of 
$\langle A_i, A_{g+i} A_j \rangle,$ 
$\langle A_i, A_jA_{g+i} \rangle,$  
$\langle A_iA_j, A_{g+i}\rangle,$ and 
$\langle A_jA_i, A_{g+i} \rangle$  
are $\Spin$, then
the argument produced in the previous lemma applies.

Suppose that $\langle A_i, A_{g+i} A_j \rangle$ is $\Pin$
and $\langle A_i A_j, A_{g+i} \rangle$ is   
$\Spin.$
This leads to the system of equations:
$ w + x + y + z = 0,$ 
$ \sqrt 2 w + \sqrt 2 z = 0,$
$w^2 + x^2 + y^2+ z^2 = 1,$ 
and
$$\tr(A_iA_j)^2 + \tr(A_{g+i})^2 +\tr( A_i A_jA_{g+i})^2 - 
\tr(A_iA_j)\tr(A_{g+i})\tr( A_i A_jA_{g+i})-2 = 2.$$
This system has non-real solutions.
A similar argument holds in the event that
$\langle A_i, A_{g+i} A_j \rangle$ is $\Spin$
and $\langle A_i A_j, A_{g+i} \rangle$ is $\Pin.$

We indicate how to proceed in the case where 
$\langle A_{i}, A_{g+i} A_j \rangle$ is $\Spin$
and $\langle A_i A_j, A_{g+i} \rangle$ is $C$.
Since $\langle A_i, A_{g+i} A_j \rangle$ is $C$,
$\tr(A_{g+i} A_j)$ and
$\tr(A_iA_{g+i}A_j)$ must take on the values $\pm \sqrt 2, 
\pm 1 $ or $0.$
The case where both are zero reduces to $\Pin$.

We explicitly work out the case 
$\tr(A_{i} A_j)= \sqrt 2,$ $\tr(A_iA_jA_{g+i})=0,$ and
$\langle A_{i}, A_{g+i} A_j \rangle$ is $\Spin.$
This leads to the equations
$$\sqrt 2 w - \sqrt 2 x = \sqrt 2$$
$$\sqrt 2 w + \sqrt 2 y=0$$
$$w^2+x^2+y^2+z^2=1$$
and
$$\tr(A_{g+i}A_j)^2 + \tr(A_i A_{g+i}A_j)^2 -\sqrt 2\tr(A_{g+i}A_j)\tr(A_iA_{g+i}A_j)= 2 .$$
This system has complex and real solutions, with
real solution $(w,x,y,z)=$ $(\frac 1 2, -\frac 1 2, -\frac 1 2, \frac 1 2),$
which implies that  $A_j \in \langle A_{i}, A_{g+i} \rangle,$ a contradiction.

Similar arguments hold in the remaining 23 cases obtained by assigning the values
$\{0, \pm 1, \pm \sqrt 2 \}$ to $\tr(A_{i} A_j)$ and $\tr(A_iA_jA_{g+i}),$
as well as the 24 cases where
$\langle A_i, A_{g+i} A_j \rangle$ is $C$
and $\langle A_i A_j, A_{g+i} \rangle$ is $\Spin$.
All cases yield systems with only non-real solutions, or
matrices $A_j \in \langle A_i, A_{g+i} \rangle.$

This leaves cases where
$\langle A_i, A_{g+i} A_j \rangle$ and $\langle A_i A_j, A_{g+i} \rangle$ are
either $C$ or $\Pin.$ In these cases, 
$\tr(A_{g+i} A_j)$, $\tr(A_{i} A_j)$, and
$\tr(A_iA_{g+i} A_j)$ must take on the values $\pm \sqrt 2,$ $\pm 1$ or $0.$


We will work out the case in which $\tr(A_{i} A_j)= \sqrt 2,$
$\tr(A_{g+i} A_j) = 1,$ and $\tr(A_iA_{g+i} A_j) = 1.$ 
The solutions to the associated system are:
$$w = \frac 3 4 \pm \frac {\sqrt{-11+8 \sqrt 2}} 4,$$
$$x = -\frac 1 4 \pm \frac {\sqrt{-11+8 \sqrt 2}} 4,$$
$$y = -\frac 1 2 \sqrt 2 + \frac 5 4 \mp \frac {\sqrt{-11+8 \sqrt 2}} 4,$$
$$z = \frac 1 2 \sqrt 2 - \frac 3 4 \mp \frac {\sqrt{-11+8 \sqrt 2}} 4,$$
with signs taken together.

However, the trace of $A_iA_{g+i}A_iA_j$ is 
$\frac 3 2 \pm \frac {\sqrt{-11+ 8 \sqrt 2}} 2.$
So the handle $\langle A_i, A_{g+i}A_iA_j \rangle$
is neither $C$ nor $\Pin$.

The remaining 123 cases are obtained by varying the
values of  $\tr(A_{i} A_j),$ $\tr(A_iA_{g+i}A_j),$ and
$\tr(A_{g+i} A_j)$ in $\{0, \pm 1, \pm \sqrt 2\}.$ The arguments proceed 
similarly to those exemplified above and
lead to: 
\begin{enumerate}
\item matrices $A_j$ inside
$\langle A_i, A_{g+i} \rangle$ (a contradiction),
\item a handle $\langle A_i, A_{g+i}A_iA_j \rangle$
that is neither $C$ nor $\Pin,$
\item one or more of $w, x, y, z$ is not real (a contradiction).
\end{enumerate}
\end{proof}

\begin{lemma}
There exists a handle $(A_i,A_{g+i})$ such that $\rho\vert_{\langle A_i, A_{g+i} \rangle}$
is neither $D$ nor $C$ nor $\Pin$.
\end{lemma}
\begin{proof}
Assume that $\langle A_i, A_{g+i} \rangle$ is $D.$
By applying Dehn twists, by the fact that $-I$ is in
the center of $\SU,$ and by the previous lemmas,
we may assume one of the following 3 cases:
\begin{enumerate}
\item $\tr(A_i) = -2s,$ $\tr(A_{g+i}) = 2s,$ and $ \tr(A_iA_{g+i})= 2s.$
\item $\tr(A_i) = 2r,$ $\tr(A_{g+i}) = -2r,$ and $ \tr(A_iA_{g+i})= -2r.$
\item $\tr(A_i) = 1,$ $\tr(A_{g+i}) = 2r,$ and $ \tr(A_iA_{g+i})= -1.$
\end{enumerate}

Case 1:

By the uniqueness of coordinates, let
$$A_i = \left[
\begin{array}{cc}
-s + \frac 1 2 i  & ir\\
ir & -s - \frac 1 2 i\\
\end{array}
\right],\ \ 
A_{g+i} = \left[
\begin{array}{cc}
s  & -r + \frac 1 2 i\\
r + \frac 1 2 i  & s\\
\end{array}
\right]. \ \ 
$$

The group $\langle A_i, A_{g+i} \rangle$
consists of: the 24 tetrahedral matrices (see above lemma),
the 8 matrices 
$$ \pm \left[
\begin{array}{cc}
s  \pm r i  &  \pm \frac 1 2 \\
\mp \frac 1 2  & s \mp r i\\
\end{array}
\right],\ \ 
$$
the 8 matrices 
$$ \pm \left[
\begin{array}{cc}
s  \pm \frac 1 2i  &  \pm r i \\
\pm r i   & s \mp \frac 1 2 i\\
\end{array}
\right],\ \ 
$$
the  8 matrices
$$ \pm \left[
\begin{array}{cc}
s   &  \pm r \pm \frac 1 2i \\
\mp r \pm \frac 1 2i  & s \\
\end{array}
\right],\ \ 
$$
(signs in the second row determined by those in the first),
and the 72 matrices obtained by multiplying those listed above by:
$$ \left[
\begin{array}{cc}
i  & 0\\
0& -i\\
\end{array}
\right],\ \ 
\left[
\begin{array}{cc}
0  & 1 \\
-1 & 0\\
\end{array}
\right],\ \ \text{ and }  
\left[
\begin{array}{cc}
0 & i\\
i & 0\\
\end{array}
\right].\ \ 
$$ 

Let $$A_j =  \left[
\begin{array}{cc}
w + ix  & y + iz\\
-y+ iz &  w - ix\\
\end{array}
\right],\ \ 
$$
with $w^2+ x^2 + y^2 + z^2 = 1$ be such that $A_j$ is
not contained in the isomorphic copy of $D$ 
generated by $\langle A_i, A_{g+i} \rangle.$
Note that
$\tr(A_iA_j) = -2s w + x - 2r z,$
$\tr(A_iA_jA_{g+i}) = - w + 2rx + 2sz,$
$\tr(A_{g+i}A_j) = 2sw + 2ry - z,$ and
$\tr(A_iA_{g+i}A_j) = -w-x-y-z.$

If each of 
$\langle A_i, A_{g+i} A_j \rangle,$ 
$\langle A_i, A_jA_{g+i} \rangle,$  
$\langle A_iA_j, A_{g+i}\rangle,$ and 
$\langle A_jA_i, A_{g+i} \rangle$  
are $\Spin$, then
the argument produced in Lemma \ref{lem:tetra} applies.

Suppose that $\langle A_i, A_{g+i} A_j \rangle$ is $\Spin$
and $\langle A_i A_j, A_{g+i} \rangle$ is $\Pin.$
This leads to the equations
$$ -2s w + x -2r z = 0,$$
$$- w + 2rx + 2sz = 0,$$
$$w^2+x^2+y^2+z^2=1,$$ and
$$\tr(A_{g+i}A_j)^2 + \tr(A_{i})^2 +\tr( A_i A_{g+i}A_j)^2 - 
\tr(A_{g+i}A_j )\tr(A_{i})\tr( A_i A_{g+i}A_j)-2 = 2.$$
The solutions are
$(w,x,y,z)= \pm (\frac 1 2, -r, 0,s)$. However, $A_j$ multiplied on the
right by
$$ \left[
\begin{array}{cc}
i  & 0\\
0& -i\\
\end{array}
\right]
$$ 
is one of the matrices listed in $\langle A_i, A_{g+i} \rangle,$
so $A_j \in \langle A_i, A_{g+i} \rangle.$

A similar argument holds in the case where
$\langle A_i, A_{g+i} A_j \rangle$ 
is $\Pin$
and $\langle A_i A_j, A_{g+i} \rangle$ is   $\Spin.$

We indicate how to proceed in the cases where 
$\langle A_{i}, A_{g+i} A_j \rangle$ is $\Spin$
and $\langle A_i A_j, A_{g+i} \rangle$ is  $D$.
Since $\langle A_{g+i}, A_{i} A_j \rangle$ is in $D$,
$\tr(A_{i} A_j)$ and
$\tr(A_iA_jA_{g+i})$ must take on the values $\pm 2s,$ $\pm 2r,$ 
$\pm 1 $ or $0.$
The case where both are zero reduces to a $\Pin$-$\Spin$ case.

We will explicitly work out the case 
$\tr(A_{i} A_j)= 2r$ and $\tr(A_iA_jA_{g+i})=-2s.$
This leads to the system:
$$-2sw + x -2rz = 2r,$$ 
$$ -w +2sz +2rx = -2s,$$
$$ w^2 +x^2 +y^2+z^2 = 1,$$ and
$$\tr(A_{g+i}A_j)^2 + \tr(A_i A_{g+i}A_j)^2 +\tr(A_i)^2 -\tr(A_i)\tr(A_{g+i}A_j)\tr(A_iA_{g+i}A_j)= 4.$$
These equations have no real solutions.

Similar arguments hold in the other 47 cases and the 48 cases where
$\langle A_i, A_{g+i} A_j \rangle$ is $D$
and $\langle A_i A_j, A_{g+i} \rangle$ is $\Spin$.
Note that some of these cases yield matrices $A_j$ that 
are in $\langle A_i, A_{g+i} \rangle.$

This leaves cases where
$\langle A_i, A_{g+i} A_j \rangle$ and $\langle A_i A_j, A_{g+i} \rangle$ are
either $D$ or $\Pin.$ In these cases, 
$\tr(A_{g+i} A_j)$, $\tr(A_{i} A_j)$, and
$\tr(A_iA_{g+i} A_j)$ must take on the values $\pm 1, \pm 2s, \pm 2r,$ or $0.$

First, we explicitly work
out a special $D$-$\Pin$ case in which $\tr(A_{i} A_j),$ $\tr(A_iA_jA_{g+i}),$ and
$\tr(A_{g+i} A_j)$ are all zero. The solutions are
$w = \mp r,$ 
$x = 0,$ 
$y = \pm \frac 1 2,$ and
$z = \pm s$ (signs taken together). 
However, the matrix $A_j$ is in $\langle A_i, A_{g+i} \rangle,$ a contradiction.

We also explicitly work out the case $\tr(A_{i} A_j)= 2r,$
$\tr(A_{g+i} A_j) = -2s,$ and $\tr(A_iA_{g+i} A_j) = 0.$ 
The three resulting equations yield $x=1.$
Thus, since $w^2 + x^2 + y^2 +z^2 = 1,$
we see that $w=y=z=0$ which violates $\tr(A_{g+i} A_j) = -2s.$

The remaining 341 cases are obtained by varying the
values of  $\tr(A_{i} A_j),$ $\tr(A_iA_{g+i}A_j),$ and
$\tr(A_{g+i} A_j)$ in $\{0, \pm 1,\pm 2r, \pm 2s\}.$ 
The arguments proceed similarly
to those exemplified above and
lead to: 
\begin{enumerate}
\item matrices $A_j$ inside
$\langle A_i, A_{g+i} \rangle$ (a contradiction),
\item a handle $\langle A_i, A_{g+i}A_iA_j \rangle$
that is not $\Pin$ nor $C$ nor $D$,
\item one or more of $w, x, y, z$ is not real (a contradiction).
\end{enumerate}

Case 2:

By the uniqueness of coordinates, let
$$A_i = \left[
\begin{array}{cc}
r -si & -\frac 1 2i\\
-\frac 1 2i & r +si\\
\end{array}
\right],\ \ 
A_{g+i} = \left[
\begin{array}{cc}
-r + \frac 1 2 i  & -s\\
s  & -r - \frac 1 2 i\\
\end{array}
\right]. \ \ 
$$

Then $\langle A_i, A_{g+i} \rangle$ yields
the 120 matrices in case 1. A study of each of the
$\Pin$-$\Spin,$ $D$-$\Spin,$ and 
$D$-$D$ cases follows as in the previous case.
As usual, each case leads to:
\begin{enumerate}
\item matrices $A_j$ inside
$\langle A_i, A_{g+i} \rangle$ (a contradiction),
\item a handle $\langle A_i, A_{g+i}A_iA_j \rangle$
that is not $\Pin$ nor $C$ nor $D$,
\item one or more of $w, x, y, z$ is not real (a contradiction).
\end{enumerate}

Case 3:

By the uniqueness of coordinates, let
$$A_i = \left[
\begin{array}{cc}
\frac 1 2 + \frac 1 2i & -\frac 1 2 - \frac 1 2i\\
\frac 1 2 - \frac 1 2i & \frac 1 2 - \frac 1 2i \\
\end{array}
\right],\ \ 
A_{g+i} = \left[
\begin{array}{cc}
r + s i  & \frac 1 2 i\\
\frac 1 2 i  & r - s i\\
\end{array}
\right]. \ \ 
$$

Then $\langle A_i, A_{g+i} \rangle$ yields
the 120 matrices in case 1. A study of each of the
$\Pin$-$\Spin,$ $D$-$\Spin,$ and 
$D$-$D$ cases follows as in case 1.
As usual, each case leads to:
\begin{enumerate}
\item matrices $A_j$ inside
$\langle A_i, A_{g+i} \rangle$ (a contradiction),
\item a handle $\langle A_i, A_{g+i}A_iA_j \rangle$
that is not $\Pin$ nor $C$ nor $D$,
\item one or more of $w, x, y, z$ is not real (a contradiction).
\end{enumerate}
\end{proof}


\begin{thebibliography}{99}
\bibitem{Be1} Benedetto, R. L., and Goldman, W. M., The Topology of the
relative character varieties of a quadruply-punctured sphere
{\em Experimental Mathematics}, 8 (1) (1999), 85-104.

\bibitem{Co1} Conway, J. H. and Jones, A. J.,  Trigonometric diophantine
equations (On vanishing sums of roots of unity) {\em Acta Arithmetica},
Vol. XXX, (1976), 229-240.

\bibitem{Ga1} Gallo, D., Kapovich, M., Marden, A.,
The Monodromy groups of schwarzian equations on closed Riemann surfaces,
{\em Preprint}.

\bibitem{Go1} Goldman, W. M.,  Ergodic theory on moduli spaces
{\em Ann. of Math.}, Vol. 146, (1997), 475-507.

\bibitem{Go2} Goldman, W. M.,  The symplectic nature of fundamental groups of 
surfaces,
{\em Adv. Math.}, Vol. 54, (1984), 200-225.
 
\bibitem{Ma1} Horowitz, R.,  Characters of free groups represented in the
two-dimensional linear group.,
{\em Comm. Pure Appl. Math.}, Vol. 25, (1972), 635-649.

\bibitem{Ma2} Magnus, W., Rings of Fricke characters and automorphism groups of free groups,
{\em Mathematische Zeitschrift}, Vol. 170, (1980), 91-103.

 
 
\bibitem{Pr1} Previte, J. P.,  and Xia, E. Z., Topological dynamics on moduli
spaces I, to appear in {\em Pacific Journal of Mathematics}.

\bibitem{Se1} Seshadri, C., Fibr\'{e}s vectoriels sur les courbes algebriques,
{\em Ast\'{e}risque }, 96 (1982).

\end{thebibliography}
\end{document}